\documentclass{amsart}

\usepackage[utf8]{inputenc}
\usepackage{amsthm, amssymb}
\usepackage[colorlinks]{hyperref}
\usepackage[colorinlistoftodos]{todonotes}
\usepackage{pgf}
\usepackage{tikz}
\usepackage{tikz-cd}
\usetikzlibrary{positioning,shapes,shadows,arrows}

\usepackage{ytableau}

\newtheorem{prop}{Proposition}
\newtheorem{thm}[prop]{Theorem}
\newtheorem{lem}[prop]{Lemma}
\newtheorem{lemma}[prop]{Lemma}
\newtheorem{cor}[prop]{Corollary}
\newtheorem{conj}[prop]{Conjecture}
\theoremstyle{remark}
\newtheorem{ex}[prop]{Example}
\newtheorem{rem}[prop]{Remark}
\newtheorem{defn}[prop]{Definition}
\numberwithin{prop}{section}
\numberwithin{equation}{section}

\newcommand{\CR}{\mathrm{CR}}
\newcommand{\CSS}{\mathrm{CSS}}
\newcommand{\Dec}{\mathrm{Dec}}
\newcommand{\df}{:=}
\newcommand{\extra}{\mathrm{ex}}
\newcommand{\fG}{\mathfrak{G}}
\newcommand{\fGb}{\fG^{(\beta)}}
\newcommand{\fL}{\mathfrak{L}}
\newcommand{\fLb}{\fL^{(\beta)}}
\newcommand{\fmin}{f_{\mathrm{min}}}
\newcommand{\fmax}{f_{\mathrm{max}}}
\newcommand{\fS}{\mathfrak{S}}
\newcommand{\Gb}{G^{(\beta)}}
\newcommand{\HD}{\xrightarrow{H}}

\newcommand{\HeckeMonoid}{\mathcal{H}}
\newcommand{\Hequiv}{\equiv_H}

\newcommand{\id}{\mathrm{id}}
\newcommand{\Inc}{\mathrm{Inc}}
\newcommand{\Insert}{\mathrm{Insert}}
\newcommand{\JNW}{J^\nwarrow}
\newcommand{\JSE}{J^\searrow}
\newcommand{\Kequiv}{\equiv_K}

\newcommand{\key}{\mathrm{key}}
\newcommand{\la}{\lambda}

\newcommand{\pab}{\partial^{(\beta)}}
\newcommand{\pib}{\pi^{(\beta)}}
\newcommand{\rev}{\mathrm{rev}}
\newcommand{\RevInsert}{\mathrm{RevInsert}}
\newcommand{\RSVT}{\mathrm{RSVT}}
\newcommand{\RSSYT}{\mathrm{RSSYT}}

\newcommand{\switch}{\mathrm{switch}}
\newcommand{\VNW}{V^\nwarrow}
\newcommand{\VSE}{V^\searrow}
\newcommand{\word}{\mathrm{word}}
\newcommand{\wt}{\mathrm{wt}}
\newcommand{\Z}{\mathbb{Z}}
\newcommand{\C}{\mathcal{C}}
\newcommand{\Ct}{\tilde{\mathcal{C}}}
\newcommand{\T}{\mathcal{T}}
\newcommand{\Tt}{\tilde{\mathcal{T}}}

\author{Mark Shimozono}
\address{Department of Mathematics \\
460 McBryde Hall, Virginia Tech\\
 255 Stanger St. \\
Blacksburg, VA, 24601, USA }
\email{mshimo@math.vt.edu}

\author{Tianyi Yu}
\address{Department of Mathematics\\
University of California, San Diego\\
La Jolla, CA, 92093, USA}
\email{tiy059@ucsd.edu}

\title{Grothendieck-to-Lascoux expansions}
\date{May 2021}

\begin{document}

\begin{abstract} 
We establish the conjecture of Reiner and Yong for an explicit combinatorial formula for the expansion of a Grothendieck polynomial into the basis of Lascoux polynomials. This expansion is a subtle refinement of its symmetric function version due to Buch, Kresch, Shimozono, Tamvakis, and Yong, which gives the expansion of stable Grothendieck polynomials indexed by permutations into Grassmannian stable Grothendieck polynomials. Our expansion is the $K$-theoretic analogue of that of a Schubert polynomial into Demazure characters, whose symmetric analogue is the expansion of a Stanley symmetric function into Schur functions. Our expansions extend to flagged Grothendieck polynomials. 
\end{abstract}

\maketitle

\section{Introduction}
The Grothendieck polynomials $\fG_w$ of Lascoux and Sch\"utzenberger \cite{LS:Groth} are explicit polynomial representatives of the $K$-classes of structure sheaves of Schubert varieties in flag varieties. Reiner and Yong \cite{ReY} conjectured an explicit combinatorial expansion of Grothendieck polynomials into the basis of Lascoux polynomials $\fL_\alpha$ \cite{Las}.
Our first main theorem (Theorem \ref{T:first main theorem}) 
gives a new\footnote{meaning, ``not stated explicitly in the literature". See Remark \ref{R:Lascoux}.} combinatorial formula for the Lascoux polynomials. This is used to prove our
second main theorem (Theorem \ref{T:second main theorem}) which
establishes the Reiner-Yong conjecture.

\subsection{Various expansions}
The Grothendieck-to-Lascoux expansion fits into a family of related expansions. The polynomials to be expanded are the cohomological and K-theoretic Schubert bases given by the Schubert polynomials $\fS_w$ and the Grothendieck polynomials $\fG_w$ respectively, and their symmetrized or stable versions, known as the Stanley symmetric functions $F_w$ and Grothendieck symmetric functions (also known as stable Grothendieck polynomials) $G_w$. These are respectively expanded into type A Demazure characters (also called key polynomials) $\kappa_\alpha$, Lascoux polynomials $\fL_\alpha$, Schur functions $s_\la$, and Grassmannian Grothendieck symmetric functions $G_\la$.

\[
\begin{matrix}
\begin{tikzcd}
\fS_w \arrow[rr,"\mathrm{symmetrize}"] \arrow[d,"\mathrm{expand}"] \arrow[d,swap,"(a)"]  && F_w \arrow[d,swap,"\mathrm{expand}"] \arrow[d,"(b)"] \\
\kappa_\alpha\arrow[rr,swap,"\mathrm{symmetrize}"] && s_\la
\end{tikzcd} \\ \\
\text{cohomology}
\end{matrix}
\qquad
\begin{matrix}
\begin{tikzcd}
\fG_w \arrow[rr,"\mathrm{symmetrize}"] \arrow[d,"\mathrm{expand}"]\arrow[d,swap,"(c)"] && G_w \arrow[d,swap,"\mathrm{expand}"] \arrow[d,"(d)"]\\
\fL_\alpha\arrow[rr,swap,"\mathrm{symmetrize}"] && G_\la
\end{tikzcd} \\ \\
\text{$K$-theory}
\end{matrix}
\]
Using the formalism of connective $K$-theory (equivalently, introducing a harmless grading parameter $\beta$ into the Grothendieck polynomial), as we do in this article, all expansions specialize to their $K$-theoretic or cohomological counterparts by setting $\beta$ to $-1$ or $0$ respectively. 

In chronological order, expansion (b) was established by Edelman and Greene \cite{EG} via a Schensted-type insertion algorithm for reduced words. The expansion (a) was found by Lascoux and Sch\"utzenberger and proved in \cite{RS}. Expansion (d) was established by Buch, Kresch, Shimozono, Tamvakis, and Yong \cite{BKSTY} via Hecke insertion, which takes Hecke words as input. Expansion (c) is the topic of this article.

The expansion coefficients have geometric significance.
The Stanley-to-Schur coefficients of the expansion (b) coincide with large rank affine Stanley to affine Schur coefficients \cite[Prop. 9.17]{LLS}, which in turn coincide with Gromov-Witten invariants for the flag variety via Peterson's Quantum Equals Affine Theorem \cite{Pet} \cite{LS:QHGr} \cite[Part 3, \S 10]{LLMSSZ}. Specializing $w$ to a Zelevinsky permutation, 
(a) and (b) give the expansion of cohomology classes of equioriented type A quiver loci \cite[Theorem 7.14]{KMS}, the latter being shown by Buch and Fulton \cite{BF} to specialize to virtually all known variants of type A Schubert polynomials.
Expansions (c) and (d) give analogous expansions in K-theory \cite{B} \cite{Mi}.

The nonsymmetric expansions are subtle refinements of their symmetric counterparts. In the symmetric expansions there is a set of tableaux in which each tableau $T$ in the set, gives a copy of $s_\la$ or $G_\la$ where $\la$ is the shape of $T$. There is a corresponding term $\kappa_\alpha$ or $\fL_\alpha$ in the nonsymmetric expansion, but an additional datum must be supplied: a composition or extremal weight $\alpha$ in the symmetric group orbit of $\la$; see \eqref{E:L sym} through \eqref{E:S sym}. Such constructions assigning a composition to a tableau go by the general name of \emph{key}. In the crystal graph of semistandard Young tableaux of shape $\la$, the left and right keys of the tableau $T$ of shape $\la$ are given by the final and initial directions of the corresponding Littelmann path whose highest weight vector is the directed line segment from the origin to $\la$. The initial direction indicates the smallest Demazure crystal containing the given tableau.

\subsection{Grothendieck and Lascoux polynomials}
The group $S_+=\bigcup_{n\ge1} S_n$ acts on $R=\Z[\beta][x_1,x_2,\dotsc]$ by permuting the variables: for $i\ge1$ let $s_i$ exchange $x_i$ and $x_{i+1}$. We define the following operators on $R$, where an element $f\in R$ (or its fraction field) denotes the operator of left multiplication by $f$.
\begin{align}
\partial_i &= (x_i-x_{i+1})^{-1}(1 - s_i) \\
\pi_i &= \partial_i x_i \\
\pab_i &= \partial_i (1+\beta x_{i+1}) \\
\pib_i &= \pab_i x_i.
\end{align}
All satisfy the braid relations for $S_+$.


We have the operator identity 
\begin{equation}\label{E:partial product rule}
  \partial_i f = \partial_i(f) + s_i(f) \partial_i\qquad\text{for all $f\in R$.}
\end{equation}
The operators satisfy the quadratic relations
\begin{align}
\partial_i^2 &= 0 \\
\pi_i^2 &= \pi_i \\
(\pab_i)^2 &= - \beta \pab_i \\
(\pib_i)^2 &= \pib_i.
\end{align}
Let $w_0^{(n)}\in S_n$ be the long element and $\rho^{(n)}=(n-1,n-2,\dotsc,1,0)$. For $w\in S_n$ the $\beta$-Grothendieck polynomial is defined by \cite{LS:Groth}
\begin{align}
\fGb_w &=
\begin{cases}
x^{\rho^{(n)}} &\text{if $w=w_0^{(n)}$} \\
\pab_i \fGb_{ws_i} &\text{if $ws_i<w$.}
\end{cases}
\end{align}
Since the $\pab_i$ satisfy the braid relations, $\fGb_w$ is well-defined for $w\in S_n$. It is also well-defined for $w\in S_+$, that is, unchanged under the standard embedding $S_n\to S_{n+1}$ for all $n\ge1$.
The Schubert $\fS_w$ and Grothendieck polynomials $\fG_w$ are defined by
\begin{align}
\fS_w &= \fGb_w|_{\beta=0}\\
\fG_w &= \fGb_w|_{\beta=-1}.
\end{align}

Let $\alpha=(\alpha_1,\alpha_2,\dotsc)$ be a composition (sequence of nonnegative integers, almost all $0$). The Lascoux polynomial $\fLb_\alpha$
is defined by \cite{Las}
\begin{align}
\fLb_\alpha = \begin{cases}
x^\alpha & \text{if $\alpha$ is a partition} \\
\pib_i \fLb_{s_i\alpha} &\text{if $\alpha_i<\alpha_{i+1}$.}
\end{cases}
\end{align}
The Demazure character $\kappa_\alpha$ is defined by
\begin{align}
\kappa_\alpha = \fLb_{\alpha}|_{\beta=0}.
\end{align}

Given a composition $\alpha=(\alpha_1,\dotsc,\alpha_n)\in\Z_{\ge0}^n$ let $\alpha^+$ be the unique partition in the $S_n$-orbit of $\alpha$. For $w\in S_n$ and $w_0\in S_n$ the long element we have the symmetrizations
\begin{align}
\label{E:L sym}
\pib_{w_0}(\fL_\alpha) &= \Gb_{\alpha^+}(x_1,\dotsc,x_n) \\
\label{E:kappa sym}
\pi_{w_0}(\kappa_\alpha) &= s_{\alpha^+}(x_1,\dotsc,x_n) \\
\label{E:G sym}
\pib_{w_0}(\fG_w(x)) &= G_w(x_1,\dotsc,x_n) \\
\label{E:S sym}
\pi_{w_0}(\fS_w(x)) &= F_w(x_1,\dotsc,x_n).
\end{align}

\subsection{New tableau formula for Lascoux polynomials}
\label{SS:new Lascoux formula}
Given the definition of a certain kind of tableau which involves entries in a totally ordered set, 
we say ``reverse" to mean the same definition but with the total order reversed.
So a \emph{reverse semistandard Young tableau} (RSSYT) is a tableau in which the entries weakly 
decrease along rows from left to right and strictly decrease along columns from top to bottom.

For a partition $\la$ a \emph{reverse set-valued tableau} (RSVT) $T$ of shape $\la$ is a filling of the boxes of $\la$ by finite subsets of $\Z_{>0}$ satisfying the following. For the box $s\in \la$ let $T(s)$ be the set which occupies the box $s$ in $T$.
\begin{enumerate}
\item $\min(T(s)) \ge \max(T(t))$ if the box $t$ is immediately right of the box $s$ in $\la$.
\item $\min(T(s))>\max(T(t))$ if the box $t$ is immediately below the box $s$ in $\la$.
\end{enumerate}
This is the reverse of Buch's set-valued tableaux \cite{B:Gr}.

Given a RSVT $T$, let $L(T)$ be the RSSYT obtained from $T$ by replacing every entry $T(s)$ by its largest value $\max(T(s))$.

The \emph{weight} $\wt(T)$ of a tableau $T$ is the composition 
whose $i$-th part is the total number of times $i$ appears in $T$.

The \emph{left key} $K_-(T)$ of a RSSYT $T$ is a composition computed by the algorithm in \S \ref{SS:families and left keys}. Up to an order-reversal bijection this is equivalent to the right key of a SSYT defined by Lascoux and Sch\"utzenberger \cite{LS:keys}; our algorithm is a variant of that of Willis \cite{Willis}.  Let $|\alpha|=\sum_{i\ge1} \alpha_i$. 
Let $\RSVT_\la$ be the set of reverse set-valued tableaux of shape $\la$.
For $T\in\RSVT_\la$ let $\extra(T)=|\wt(T)|-|\la|$.

Our first main theorem is:

\begin{thm} \label{T:first main theorem}
For any composition $\alpha$
\begin{align}\label{E:Lascoux RSVT}
  \fLb_\alpha = \sum_{\substack{T\in\RSVT_{\alpha^+} \\ K_-(L(T))\le \alpha}} \beta^{\extra(T)} x^{\wt(T)}
\end{align}
Here $\le$ indicates the quotient of Bruhat order on the orbit $S_+\alpha$.
\end{thm}
Another way to state the $\le$ relation on compositions is the following.
A \emph{key tableau} (or just key) is a SSYT\footnote{Sometimes we will write a key tableau as a RSSYT depending on context.} of partition shape such that the $j$-th column, viewed as a set, contains the $(j+1)$-th for all $j$. There is a bijection
$\alpha\mapsto \key(\alpha)$ from compositions to key tableaux where
$\key(\alpha)$ is the unique SSYT of shape $\alpha^+$ and weight $\alpha$. Its $j$-th column consists of the numbers $\{i\mid \alpha_i \ge j\}$ for all $j$.

Then for compositions $\alpha$ and $\beta$, $\alpha\le\beta$ if and only if $\alpha^+=\beta^+$ and $\key(\alpha)$ is entrywise less than or equal to $\key(\beta)$.

Theorem \ref{T:first main theorem} is proved in \S \ref{S:formula for Lascoux polynomials}.

\begin{rem} \label{R:Lascoux}
There have been a number of conjectural combinatorial formulas for Lascoux polynomials, such as the K-Kohnert move rule of Ross and Yong \cite{RoY} (\cite[Footnote on p. 19]{Kir} for the general $\beta$ version),
the set-valued skyline filling formula of Monical \cite{Mon}, and a set-valued tableau (SVT) rule of Pechenik and Scrimshaw \cite{PS} which requires the fairly involved Lusztig involution on a crystal structure on SVTs in addition to an entrywise minimum and the usual right key of a SSYT. Buciumas, Scrimshaw and Weber \cite{BSW} proved the last two of these rules using 
solvable lattice models. 
In response to a previous version of this article, Travis Scrimshaw kindly informed us that Theorems \ref{T:first main theorem} and \ref{T:first main theorem restated} are implicit in \cite{BSW}: see the proof of \cite[Theorem 4.4]{BSW}. We feel it is worthwhile to state these theorems in their simplest and most explicit form.
We note that the naive nonreversed SVT analogue of the RSVT formula does not yield the Lascoux polynomial.
\end{rem}

Equation \eqref{E:Lascoux RSVT} can be restated using only RSSYTs, avoiding SVT altogether.
Let $\RSSYT_\la$ be the set of reverse semistandard tableaux of shape $\la$.

\begin{thm} \label{T:first main theorem restated}
For any composition $\alpha$
\begin{align}
\fLb_\alpha = \sum_{\substack{T\in\RSSYT_{\alpha^+} \\ K_-(T)\le \alpha}} x^{\wt(T)} \prod_{(s,k)} (1+\beta x_k)
\end{align}
where the product runs over pairs $(s,k)$ where $s$ is a box of $\alpha^+$, $k<T(s)$, and replacing the $s$-th entry of $T$ by $k$ results in a RSSYT.
\end{thm}

\begin{rem} The naive nonreversed SSYT analogue of Theorem \ref{T:first main theorem restated} also does not give the 
Lascoux polynomial.
\end{rem}

\begin{rem} The condition on the pairs $(s,k)$ should be compared with the formula for Grothendieck polynomials indexed by vexillary permutations in part 3 of the second Corollary in section 1.2 of \cite{KMY}.
\end{rem}

\subsection{Stable limit of Lascoux polynomials}
For a fixed composition $\alpha$ consider the limit
$\lim_{N\to\infty} \fLb_{(0^N,\alpha)}$ in which more zero parts are prepended to $\alpha$. In Theorem \ref{T:first main theorem} it is evident that the limit depends only on $\la=\alpha^+$ and it is given by removing the left key condition. By the definition of Lascoux polynomial,
this limit can be computed by $\pib_i$ operators on $x^\la$. One may show that 
$$\fLb_{w_0^{(n)}\la} = \Gb_{\la}(x_1,\dotsc,x_n)\qquad\text{for $\la=(\la_1,\dotsc,\la_n)$}$$
where $\Gb_\la$ is the Grassmannian Grothendieck symmetric function, which will be defined in \S \ref{SS:Grass Groth}.

We deduce that the above limit coincides with the symmetric series $\Gb_\la$ on one hand and the generating function of RSVT of shape $\la$ on the other. This is the reversed version of Buch's SVT formula for $\Gb_\la$ \cite{B:Gr}; see \eqref{E:Grass Groth tableau}.

\subsection{Fomin-Kirillov monomial formula}
\label{SS:FK}
Our point of departure is the explicit monomial expansion of Grothendieck polynomials due to Fomin and Kirillov \cite{FK}.

The \emph{$0$-Hecke monoid} $\HeckeMonoid$ is the quotient of the free monoid of words on the alphabet $\Z_{>0}$ by the relations 
\begin{align}
i (i+1) i &\Hequiv (i+1) i (i+1) \\
ii &\Hequiv i \\
ij &\Hequiv ji \qquad\text{for $|i-j| \ge 2$.}
\end{align}
$\HeckeMonoid$ acts on $S_+$ by 
\begin{align*}
  i * w = \begin{cases}
  s_i w & \text{if $s_i w>w$} \\
  w &\text{if $s_iw<w$.}
  \end{cases}
\end{align*}
Given a word $u\in \HeckeMonoid$ define its \emph{associated permutation} by
$u * \id\in S_+$. For $w\in S_+$ let $\HeckeMonoid_w$ be the words in 
$\HeckeMonoid$ with associated permutation $w$. The subsets $\HeckeMonoid_w\subset \HeckeMonoid$ are the $\Hequiv$-equivalence classes.

\begin{lem} \label{L:reverse Hecke word} $u\in \HeckeMonoid_w$ if and only if $\rev(u)\in\HeckeMonoid_{w^{-1}}$.
\end{lem}

For $a\in\HeckeMonoid_w$ let $\extra(a) = \text{length}(a) - \ell(w)$, the \emph{excess} of the length of $a$ above the minimum possible, the Coxeter length $\ell(w)$ of $w$.

The following is merely the definition in \cite{BJS} but with both words reversed, which is better suited to our 
use of decreasing tableaux.
\begin{defn} \label{D:compatible pair} \cite{BJS}
A pair of words $(a,i)$
is \emph{compatible} if
they satisfy
\begin{enumerate}
\item $a,i$ are words of positive
numbers with the same length.
\item $i$ is weakly decreasing
\item $i_j = i_{j+1}$ implies 
$a_j < a_{j+1}$.
\end{enumerate}
A compatible pair $(a,i)$ is \emph{bounded} if
\begin{align}
\label{E:bounded}
    i_j \le a_j\qquad\text{for all $j$.}
\end{align}
\end{defn}
Let $\C$ be the set of all compatible pairs, $\C^b$ those that are bounded, $\C_w$ the pairs $(a,i)\in \C$ such that $a\in \HeckeMonoid_w$, and $\C_w^b = \C^b \cap \C_w$.

The following monomial expansion of $\beta$-Grothendieck polynomials is due to Fomin and Kirillov \cite{FK}.
\begin{align}\label{E:Groth compatible pairs}
\fGb_w = \sum_{(a,i)\in \C_{w^{-1}}^b} \beta^{\extra(a)} x^{\wt(i)}
\end{align}
When $\beta=0$ this is the Billey-Jockusch-Stanley formula for Schubert polynomials \cite{BJS}.

For $w\in S_n$ and a positive integer $N$ let
$1^N\times w$ be the permutation of $S_{n+N}$ obtained by adding $N$ fixed points before $w$. The $\beta$-Grothendieck symmetric function is defined by
\begin{align}
\Gb_w = \lim_{N\to\infty} \fG_{1^N\times w}
\end{align}
It lives in a completion of the ring of symmetric functions over $\Z[\beta]$.
The Stanley and Grothendieck symmetric functions are defined by
\begin{align}
F_w &= \Gb_w|_{\beta=0} \\
G_w &= \Gb_w|_{\beta=-1}.
\end{align}
It follows from \eqref{E:Groth compatible pairs} and the definitions that
\begin{align}\label{E:GSF compatible pairs}
\Gb_w = \sum_{(a,i)\in \C_{w^{-1}}} \beta^{\extra(a)} x^{\wt(i)}.
\end{align}
Comparing with \eqref{E:Groth compatible pairs} just the boundedness condition has been dropped.

\subsection{$\Gb_w$ to $\Gb_\la$ via Hecke insertion: restriction of compatible pairs according to $w$}
\label{SS:Grass Groth}
The code $c(w)$ of a permutation is the sequence $(c_1,c_2,\dotsc)$ such that 
$c_i = | \{j \mid \text{$1 \le j < w(i)$ and $w^{-1}(j)>i$} \}|$.
For a partition $\la=(\la_1,\la_2,\dotsc,\la_k)$ the \emph{Grassmannian Grothendieck symmetric function} $\Gb_\la$ is by definition equal to $\Gb_w$ where $w$ is the permutation with code $(\la_k,\dotsc,\la_2,\la_1,0,0,\dotsc)$.

Buch \cite{Bu} showed that the $\Z[\beta]$-span of the $\Gb_w$ for $w\in S_+$, has basis given by the $\Gb_\la$ and proved the increasing version of the following:
\begin{align}\label{E:Grass Groth tableau}
\Gb_\la = \sum_{T\in\RSVT_\la} \beta^{\extra(T)} x^{\wt(T)}
\end{align}
For $\beta=0$ this becomes the RSSYT formula for the Schur function $s_\la$.

To find the coefficients of the $\Gb_w$ to $\Gb_\la$ expansion, the \emph{Hecke insertion} algorithm was developed in \cite{BKSTY} in the language of \emph{increasing tableaux}, which strictly increase along rows from left to right and 
strictly increase along columns from top to bottom. In \S \ref{SS:Hecke insertion} we recall these definitions but use the variant of Hecke insertion for \emph{decreasing tableaux}, which are the ``reverse" of increasing: they strictly decrease along rows from left to right and strictly decrease along columns from top to bottom.

It was not explicitly stated in \cite{BKSTY} but all ingredients are there to define a Hecke Robinson-Schensted-Knuth (RSK) bijection called Insert (and its inverse bijection RevInsert)
\tikzstyle{start}=[to path={(\tikztostart.#1) -- (\tikztotarget)}]
\[
\begin{tikzcd}[every arrow/.append style={shift left}]
    \C \arrow[rr,"\textrm{Insert}"] && \bigsqcup_\la (\Dec_\la \times \RSVT_\la)=:\T \arrow[ll,"\textrm{RevInsert}"] \\
    (a,i) \arrow[rr] && (P,Q) \arrow[ll]
\end{tikzcd}
\]
where $\Dec_\la$ is the set of decreasing tableaux of shape $\la$. Note that the set $\T$ as defined by the above diagram, consists of pairs $(P,Q)$ of tableaux of the same partition shape with $P$ decreasing and $Q$ reverse set-valued.

Let $(P,Q)=\Insert(a,i)$. 
By Proposition \ref{P:Hecke insertion Kequiv and weight} the bijection satisfies
\begin{align}
\label{E:insert K equiv}
\rev(a) &\Kequiv P \\
\label{E:insert wt preserving}
\wt(Q)&=\wt(i)
\end{align}
where the relation $\Kequiv$ is defined in \S \ref{SS:right key}. Note that $\Kequiv$ refines the relation $\Hequiv$ of \S \ref{SS:FK}. 
By Lemma \ref{L:reverse Hecke word} and \eqref{E:insert K equiv} the  bijection $\Insert$ restricts to a bijection 
\begin{align}\label{E:w bijection}
\C_{w^{-1}}\leftrightarrow \bigsqcup_\la \left(\Dec_\la^w \times \RSVT_\la\right) =: \T_w
\end{align}
where $\Dec_\la^w = \{ T\in \Dec_\la\mid \word(T)\in \HeckeMonoid_w\}$ and $\word(T)$ is defined in \S \ref{SS:background}; see \cite{BKSTY} for the increasing tableau version.
Taking the generating function of both sides we obtain
\begin{align}
  \Gb_w = \sum_\la |\Dec_\la^w| \Gb_\la.
\end{align}

\subsection{$\fGb_w$ to $\fLb_\alpha$ by Hecke insertion and keys: restriction to bounded compatible pairs}
\label{SS:bounded bijection}
We give an algorithm to compute the right key $K_+(T)$ of a RSSYT $T$ in Definition \ref{D:right key of RSSYT}. It is essentially equivalent to that of Willis \cite{Willis} in the context of usual RSK and Knuth equivalence. However we use it for a new purpose, applying it to decreasing tableaux in the context of Hecke insertion and the K-Knuth equivalence $\Kequiv$.
The \emph{K-jeu-de-taquin} (Kjdt) of Thomas and Yong \cite{TY} may be used (following \cite{ReY} for increasing tableaux) to give another definition of right key of decreasing tableau, 
and two definitions of right key are shown to coincide (Proposition \ref{P:Kjdt right key}). 
For our amusement we give a third formulation of right key for a decreasing tableau in \S \ref{S:alternative right key} using the ``transpose" of Hecke reverse insertion; in the context of semistandard tableaux and Knuth equivalence, this formulation is analogous to an original definition of Lascoux and Sch\"utzenberger \cite{LS:keys}.

Let $\T^b$ be the subset of pairs $(P,Q)\in \T$ such that $K_+(P) \ge K_-(L(Q))$.
With $\T_w$ defined as in \eqref{E:w bijection} let $\T_w^b = \T^b \cap \T_w$.

In \S \ref{S:compatible} we show the following (Theorem \ref{restricting_Insert}):

\begin{thm} \label{T:compatible bounded}
$\Insert$ restricts to a bijection $\C^b \cong \T^b$.
\end{thm}

Intersecting with the bijection \eqref{E:w bijection}, $\Insert$ restricts to a bijection 
\begin{align}\label{E:w bounded bijection}
\C^b_{w^{-1}}\cong \T^b_w \qquad\text{for every $w\in S_+$.} 
\end{align}
Using Theorem \ref{T:first main theorem} we obtain our second main theorem, 
the Grothendieck-to-Lascoux expansion via decreasing tableaux.

\begin{thm} \label{T:second main theorem} 
\begin{align}\label{E:Groth to Lascoux}
\fGb_w = \sum_\la \sum_{P\in \Dec_\la^w} \fLb_{\wt(K_+(P))}.
\end{align}
\end{thm}

\subsection{Connecting with the Reiner-Yong conjecture}

The Reiner-Yong conjecture asserts:

\begin{conj} \cite{ReY} 
\begin{align}
\fGb_w = \sum_\la \sum_{P\in \Inc_\la^{w^{-1}}} \fLb_{\wt(K_-(P))}
\end{align}
where $\Inc_\la^w$ is similar to $\Dec_\la^w$ except that the tableaux are increasing and $K_-(P)$ is the left key construction on the increasing tableau $P$ using the Kjdt. 
\end{conj}
\begin{proof} By Propositions \ref{P:decreasing to increasing} and \ref{reverse_complement_fixes_keys} the map $T\mapsto T^\sharp$ defines a bijection $\Dec_\la\to\Inc_\la$ such that $K_+(T)=K_-(T^\sharp)$ and $\word(T^\sharp)\Kequiv \rev(\word(T))$.
Using Lemma \ref{L:reverse Hecke word} we see that $T\in \HeckeMonoid_w$ if and only if $T^\sharp\in\HeckeMonoid_{w^{-1}}$. Thus 
the bijection restricts to a bijection $\Dec_\la^w\cong \Inc_\la^{w^{-1}}$ as required.
\end{proof}

\subsection{Cohomological case}
Say that a word is \emph{reduced} if it is a reduced word for some permutation in $S_+$. Say that a tableau is reduced if its word is. Let $\C^b(0)$ denote the set of compatible pairs $(a,i)$ such that $a$ is reduced.
Let $\T^b(0)$ be the set of pairs $(P,Q)$ of tableaux of the same partition shape such that $P$ is a reduced decreasing tableau and $Q$ is a RSSYT with $K_+(P) \ge K_-(Q)$. 
By setting $\beta=0$ in Theorem \ref{T:compatible bounded} we obtain the following.

\begin{thm} \label{T:compatible bounded reduced} 
\begin{enumerate}
    \item The decreasing analogue of Lascoux and Sch\"utzenberger's right nil key of a reduced decreasing tableau $P$ (see left nil key for reduced increasing tableaux in \cite{RS}) agrees with $K_+(P)$.
    \item The restriction of $\Insert$ gives a bijection $\C^b(0)\cong \T^b(0)$. In this case $\Insert$ is computed by
Edelman-Greene column insertion into reduced decreasing tableaux, recorded by RSSYT.
\end{enumerate}

\end{thm}

This recovers the Schubert to Demazure expansion \cite{LS:Schub} \cite{RS}. 

\begin{rem}\label{R:Schubert Demazure crystal}
The Demazure crystal structure on SSYT was essentially known to Lascoux and Sch\"utzenberger \cite{LS:keys}.
Theorem \ref{T:compatible bounded reduced} clarifies the Demazure crystal structure in the Schubert expansion: it is pulled back via Edelman-Greene column insertion from the Lascoux-Sch\"utzenberger Demazure crystal structure on the semistandard $Q$-tableau. See also \cite{AS}.
\end{rem}

\subsection{Plactic variant}
\label{SS:plactic variant}

Let $\Ct$ be the set of Knuth biwords,
which are defined by changing the third condition of Definition \ref{D:compatible pair} to: $i_j=i_{j+1}$ implies $a_j \le a_{j+1}$. Let $\Tt$ be the set of pairs $(P,Q)$ of RSSYT of the same partition shape. Let $\Ct^b$ be the bounded compatible pairs in $\Ct$ and $\Tt^b$ the pairs $(P,Q)\in \Tt$ such that
$K_+(P) \ge K_-(Q)$ using the classical definition of left and right key in \cite{LS:keys} 
but adapted to RSSYT.

\begin{thm} The value-reversed version of column insertion RSK yields a bijection 
$\Ct^b \cong \Tt^b$.
\end{thm}

See Remark \ref{R:semistandard compatible bounded} for more details. This result applies to 
the representation theory of the general linear group over a non-archimedean local field \cite{GLS}.

\subsection{Flagged Grothendieck to Lascoux}
In this subsection we extend our expansion to flagged Grothendieck polynomials.

In the literature there is a definition of flagged Grothendieck polynomial whose generality extends to the 
case of 321-avoiding permutations \cite{Mat}; see \cite{KMY2} for the case of vexillary permutations.
For 321-avoiding permutations there is a monomial tableau formula and a determinantal formula. 

We use a divided difference definition of flagged Grothendieck polynomial from \cite{LLS:back stable Groth} 
which is valid for any permutation. This flagged Grothendieck polynomial has an explicit monomial expansion given in Proposition \ref{P:flagged G}.

The main result of this subsection is a Lascoux polynomial expansion of flagged Grothendiecks.

An \emph{flag} is a sequence of integers $f=(f_1,f_2,\dotsc,f_n)$ which is weakly increasing, satisfies 
$f_i \ge i$ for all $i$, and $f_n=n$. Let $\fmin=(1,2,\dotsc,n)$ and $\fmax=(n,n,\dotsc,n)$ be the minimum and maximum flags respectively. Given a flag $f$, define the permutation $\sigma_f\in S_n$ as follows. For the minimum flag $\fmin$ let $\sigma_{\fmin} = \id$. For $f\ne \fmin$ there is an index $j$ such that $f_j > j$; take the minimum such. Define $\sigma_f = s_i \sigma_{f'}$ where $i+1=f_j$ and $f'$ is obtained from $f$ by replacing the $i+1$ by $i$. The \emph{flagged Grothendieck polynomial} is defined by $\fGb_{w,f}=\pib_{\sigma_f}(\fGb_w)$.

The flagged Grothendieck polynomials have the following explicit monomial expansion.

\begin{prop}\label{P:flagged G}\cite{LLS:back stable Groth}
	\begin{align}\label{E:flagged G}
		\fGb_{w,f} &= \sum_{\substack{(a,i)\in \C_w \\ i_k \le f_{a_k}}}
		\beta^{\extra(a)} x^{\wt(i)} 
	\end{align}
\end{prop}

\begin{rem}
Note that only the bound $i_k \le a_k$ in \eqref{E:Groth compatible pairs} has been changed to 
$i_k \le f_{a_k}$. 
\end{rem}

The flagged Grothendieck polynomials interpolate between Grothendieck polynomials and 
their symmetric counterparts.

\begin{cor} \label{C:symmetrize G} For $w\in S_n$
	\begin{align}\label{E:symmetrize G}
		\pib_{w_0}(\fGb_w) = \Gb_w(x_1,\dotsc,x_n).
	\end{align}
\end{cor}
\begin{proof} We have
\begin{align*}
\pib_{w_0}(\fGb_w) 
&= \pib_{w_0}(\fGb_{w,\fmin}) \\
&= \fGb_{w,\fmax} \\
&= \Gb_w(x_1,\dotsc,x_n)
\end{align*}
where the last equality holds by the equality of \eqref{E:flagged G} with
\eqref{E:GSF compatible pairs} with $x_i$ set to $0$ for $i>n$.
\end{proof}

Define the Demazure action $\circ$ of $S_+$ on compositions by
\begin{align}
   s_i \circ \alpha = \begin{cases}
   s_i(\alpha) & \text{if $\alpha_i > \alpha_{i+1}$} \\
   \alpha & \text{otherwise.}
   \end{cases}
\end{align}

Theorem \ref{T:second main theorem} implies the following.

\begin{cor} 
\begin{align}
    \fGb_{w,f} = \sum_\lambda \sum_{P\in \Dec_\lambda^w} \fLb_{\sigma_f \circ \wt(K_+(P))}.
\end{align}
\end{cor}

\subsection*{Acknowledgements}
The authors thank Alex Yong for helpful conversations and especially for sharing with us the details of his conjecture with Vic Reiner. M. S. thanks Tomoo Matsumura for help related to flagged Grothendieck polynomials. T. Y. thanks Brendon Rhoades for helpful conversations.

\section{Background}

\subsection{Partitions, tableaux, words}
\label{SS:background}
A skew shape is \emph{normal} (resp. \emph{antinormal}) if it is empty or has a unique northwestmost (resp. southeastmost) corner. Given a skew shape $D$ let $D^*$ be a skew shape obtained by 180 degree rotation of $D$.

For a tableau $T$ let $\word(T)$ be the \emph{column-reading word} of $T$, obtained by reading the first column of $T$ from bottom to top, then the second column of $T$ from bottom to top, and so on.Thus the word of a column of an increasing (resp. decreasing) tableau is a decreasing (resp. increasing) word. We shall often make no distinction between a column, its word, and the underlying set. Similarly we may write $T$ and mean $\word(T)$.

\subsection{Thomas and Yong's $K$-theoretic jeu-de-taquin (Kjdt)}

\subsubsection{Kjdt}

We will define the Kjdt using the language of increasing tableaux as it is more natural. But the same definitions apply to decreasing tableaux in the same way; one is just labeling the chains of partitions differently.

Say that the skew shapes $D$ and $E$ are \emph{successive} if there are partitions $\nu\subset\mu\subset\la$ such that $D=\mu/\nu$ and $E=\la/\nu$. In that case let $D+E=\la/\nu$. We also say that $E$ \emph{extends} $D$.

A \emph{rook strip} is a skew shape that is both a \emph{horizontal strip} (has at most one box in each column) 
and a \emph{vertical strip} (has at most one box in each row). Say that a skew shape is \emph{thin} if it is the sum of two successive rook strips. Note that a skew shape is thin if and only if it has no $2\times 2$ subdiagram and has at most two boxes in any row and column. When decomposing a thin shape into a sum of successive rook strips the only choice is for isolated boxes, which can be in either the inner or outer rook strip. 

Let $D$ and $E$ be successive rook strips. The \emph{switch} of $(D,E)$ is the unique successive pair of rook strips $(E',D')$ such that $E'+D'=D+E$, every isolated box of $D$ is in $D'$ and every isolated box of $E$ is in $E'$.
\begin{ex}\label{X:switch} In the following the boxes of $D$ and $D'$ are filled with $\circ$ and those of 
$E$ and $E'$ are filled with $\bullet$.
\ytableausetup{aligntableaux=center}
\begin{align*}
D+E=
\ytableaushort{\none\none\none\none\bullet,\none\none\none\circ,\none\none\circ\bullet,\none\circ\bullet,\circ}
\leftrightarrow
\ytableaushort{\none\none\none\none\bullet,\none\none\none\bullet,\none\none\bullet\circ,\none\bullet\circ,\circ}=E'+D'
\end{align*}
\end{ex}
An increasing skew tableau $T$ of shape $D$ can be viewed as a sequence of successive rook strips
$D_1,D_2,\dotsc,D_n$ with $D=D_1+D_2+\dotsm+D_n$: the boxes of $D_i$ are filled with $i$.
Let $E$ be a rook strip which extends $D$. The \emph{forward Kjdt} $\JSE_E(T)$ is the increasing tableau defined by switching $E$ past all the $D_i$. More precisely let $E^{(n)} = E$ and define
$\switch(D_i,E^{(i)}) = (E^{(i-1)},D'_i)$ for all $i$ going from $n$ down to $1$. Then $\JSE_E(T)$ is the increasing tableau defined by the successive rook strips $D_1',D_2',\dotsc,D_n'$.
The \emph{vacated rook strip} $\VSE_E(T)$ is by definition $E^{(0)}$.

\begin{ex} \label{X:JSE one strip}
Here is an example of $\JSE_E(T)$ where the boxes of $E^{(i)}$ are filled with $\bullet$
and $D_i$ and $D_i'$ are filled with $i$'s. $E$ is filled with $\bullet$ in the first diagram. $T$ and $T'$ are the  skew increasing tableaux in the first and last diagrams respectively, where one ignores the $\bullet$'s.
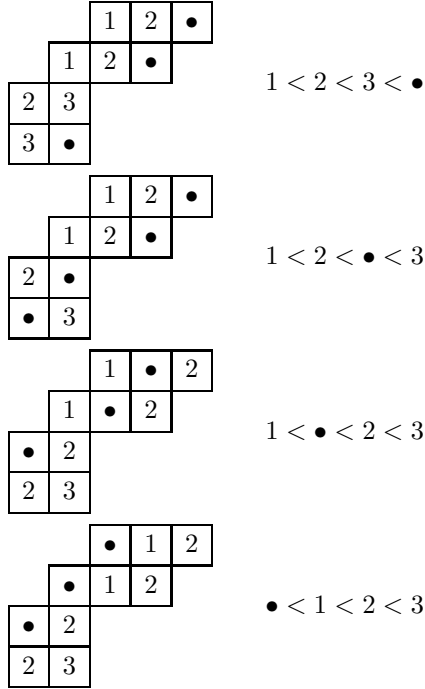
\begin{align*}
&\ytableaushort{\none\none12\bullet,\none12\bullet,23,3\bullet} \qquad 1 < 2 < 3 < \bullet \\
&\ytableaushort{\none\none12\bullet,\none12\bullet,2\bullet,\bullet3} \qquad 1 < 2 < \bullet < 3 \\
&\ytableaushort{\none\none1\bullet2,\none1\bullet2,\bullet2,23} \qquad 1 < \bullet < 2 < 3 \\
&\ytableaushort{\none\none\bullet12,\none\bullet12,\bullet2,23} \qquad \bullet < 1 < 2 < 3 \\
\end{align*}
Note that we can think of all the intermediate ``tableaux" as being increasing with respect to the total orders indicated to the right.
\end{ex}

Suppose $T$ and $U$ are increasing tableaux whose shapes are successive.
Let $T$ have rook strip decomposition $D_1+D_2+\dotsm+D_n$
and let $U$ be $E_1+E_2+\dotsm+E_m$. Then define $\JSE_U(T)$ to be the increasing tableau given by
switching $T$ past $E_1$ using $\JSE_{E_1}$ then past $E_2$ by $\JSE_{E_2}$ etc. Let $\VSE_U(T)$ be the sequence of vacated rook strips.
Both are increasing tableaux.

\begin{ex} \label{X:JSE} Using $T$ from Example \ref{X:JSE one strip} let $U$ be the increasing skew tableau filled with $a$'s and $b$'s. $T$ is first switched past the $a$'s as in Example \ref{X:JSE one strip} 
and then past the $b$'s.
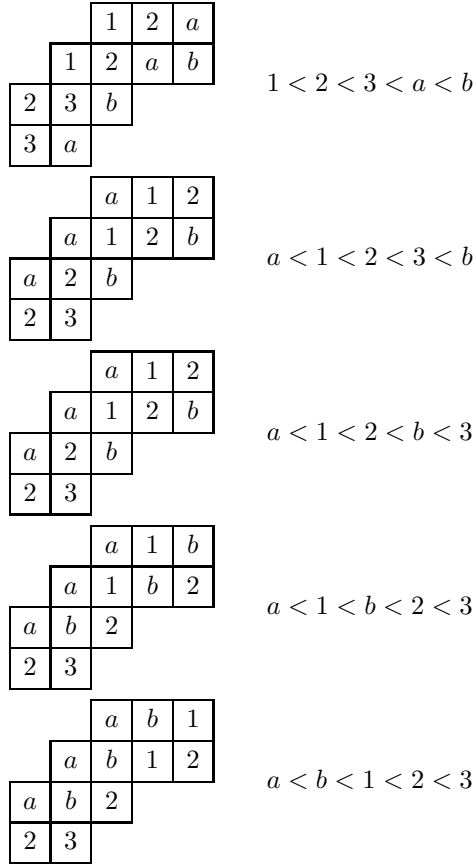
\begin{align*}
&\ytableaushort{\none\none12a,\none12ab,23b,3a} \qquad 1 < 2 < 3 < a < b \\
&\ytableaushort{\none{\none}a12,{\none}a12b,a2b,23} \qquad a < 1 < 2 < 3 < b\\
&\ytableaushort{\none{\none}a12,{\none}a12b,a2b,23} \qquad a < 1 < 2 < b < 3\\
&\ytableaushort{\none{\none}a1b,{\none}a1b2,ab2,23} \qquad a < 1 < b < 2 < 3\\
&\ytableaushort{\none{\none}ab1,{\none}ab12,ab2,23} \qquad a < b < 1 < 2 < 3\\
\end{align*}
The final resting place of the numbers is $\JSE_U(T)$ and that of the letters is $\VSE_U(T)$.
\end{ex}

Instead of switching into rook strips that are to the outside of our skew increasing tableau, we can switch into rook strips lying on the inside. Given a rook strip $D$ and increasing tableau $U$ of a shape $E$ which extends $D$, let $\JNW_D(U)$ be the increasing skew tableau obtained by switching $D_1$ with the first rook strip of $U$, then switching with the second rook strip of $U$, and so on. Let $\VNW_D(U)$ be the rook strip vacated by this process. Here $U$ moves to the northwest. Similarly if $T$ and $U$ are skew increasing tableaux whose shapes are successive, we may define $\JNW_T(U)$ by switching $U$ to the northwest, past the last strip of $T$, then the next to last, and so on. Let $\VNW_T(U)$ be the increasing tableau defined by the sequence of rook strips vacated by this process.

It is equivalent to think of the inner tableau as a sequence of switching instructions to apply to the outer, as it is to think of the outer tableau as a sequence of switching instructions to apply to the inner. Both give the same result.

\begin{thm} \cite[Theorem 3.1]{TY} 
Let $T$ and $U$ be increasing tableaux whose shapes are successive. 
Then $\JSE_U(T)=\VNW_T(U)$ and $\VSE_U(T)=\JNW_T(U)$.
\end{thm}

Let $T$ and $U$ be increasing tableaux whose shapes are successive. The \emph{infusion} of the pair $(T,U)$ is the pair $(\VSE_U(T),\JSE_U(T)) = (\JNW_T(U),\VNW_T(U))$.

\begin{rem} 
\begin{itemize}
    \item
More generally is known (see \cite{BSS} \cite{H} for jeu-de-taquin for semistandard tableaux) that infusion can be computed by any sequence of switches which ``shuffles" the alphabets of the inner and outer tableaux. In Example \ref{X:JSE} this means starting with $1<2<3<a<b$, ending with $a<b<1<2<3$ and always using a total order on these 5 values with $1<2<3$ and $a<b$.     
    \item
Since a single switch is an involution it follows that infusion is involutive.    
\end{itemize}

\end{rem}

\subsubsection{$K$-Pieri property of Kjdt}
A \emph{horizontal $K$-Pieri $t$-strip} (called a $t$-Pieri filling in \cite{TY}) is sequence of $t$ nonempty successive rook strips $D_1,D_2,\dotsc,D_t$ such that $D_1+D_2+\dotsm+D_t$ is a horizontal strip and the boxes of $D_{i+1}$ are to the right of the boxes of $D_i$ for all $1\le i< t$. This can be depicted by an increasing tableau in which the boxes of $D_i$ are filled with $i$. A \emph{vertical $K$-Pieri strip} is the transpose analogue.

\begin{ex}\label{X:K Pieri strip} A horizontal $K$-Pieri $4$-strip is pictured below as an increasing tableau.
\begin{align*}
\ytableaushort{
\none\none\none\none\none\none34,\none\none\none\none\none2,\none\none\none\none,\none\none\none2,\none\none2,12}
\end{align*}
\end{ex}

\begin{prop}\label{P:Kjdt Pieri} \cite{TY} The Kjdt operations $\JNW_T$ and $\JSE_U$ send horizontal (resp. vertical) $K$-Pieri $t$-strips, to horizontal (resp. vertical) $K$-Pieri $t$-strips.
\end{prop}

\subsubsection{Two special rectification orders}
In the infusion of the pair $(T,U)$, as computed by $\JNW_T(U)$ (resp. $\JSE_U(T))$, the tableau $T$ (resp. $U$) is called a \emph{rectification order}; the tableau $T$ (resp. $U$) is viewed as the sequence of instructions for the moving of $U$ (resp. $T$) to the northwest (resp. southeast). 
We require two special kinds of tableaux for this purpose.

Given a partition $\la$ let $\CSS(\la)$ be the \emph{column superstandard tableau}, the tableau of shape $\la$ with first column is filled with $1$ through $\la_1'$ from top to bottom, second column filled with $\la_1'+1$ through $\la_1'+\la_2'$,  and so on. 
For $\la'=(\la_1',\dotsc,\la_k')$
let $\CR(\la)$ be the \emph{column reading tableau}, the unique tableau of 
shape $\la$ built from the empty tableau by placing the numbers $1$ through $\la_k'$ at the ends of the first $\la_k'$ rows, then placing the next $\la_{k-1}'$ numbers at the ends of the first $\la_{k-1}'$ rows, and so on.

We make similar definitions of column superstandard and column reading antitableaux $\CSS(\la^*)$ and $\CR(\la^*)$, which can also be obtained respectively from $\CSS(\la)$ and $\CR(\la)$ by rotation by 180 degrees and complementing. 

We shall sometimes use decreasing versions of these special tableaux.

\begin{ex} For $\la=(4,3,2)$ we have $\la'=(3,3,2,1)$ and  
\begin{align*}
\CSS(\la) &= \ytableaushort{1479,258,36} &\qquad
\CR(\la) &= \ytableaushort{1247,358,69} \\
\CSS(\la^*) &= \ytableaushort{\none\none47,\none258,1369} &
\CR(\la^*) &= \ytableaushort{\none\none14,\none257,3689}
\end{align*}
\end{ex}

\subsubsection{$K$-rectification and anti-rectification}
In studying the Kjdt for most situations it does not matter whether a tableau is increasing or decreasing because the tableaux are specified by sequences of rook strips and only the labeling of the strips is different. For technical reasons we will use different rectification orders in the following definitions for increasing versus decreasing tableaux.

\begin{defn}
\begin{enumerate}
\item 
Given an increasing (resp. decreasing) tableau $T$ of shape $D=\la/\mu$ define its \emph{$K$-rectification} by $\JNW(T) = \JNW_S(T)$ where $S=\CR(\mu)$ (resp. $S=\CSS(\mu)$). 
\item
Let $R$ be a tight ($\ell(\la) \times \la_1$) rectangle placed around $\la$. Then $R/\la$ is an antinormal skew shape. Define the \emph{$K$-anti-rectification}\footnote{We feel this name is more descriptive as the result has antinormal shape. The original name is ``reverse K-rectification" \cite{TY}.} of the increasing (resp. decreasing) tableau $T$ by $\JSE(T) = \JSE_U(T)$ where $U=\CSS(R/\la)$ (resp. $U=\CR(R/\la)$).
\item Given a decreasing tableau $T$ of partition shape define $T^\sharp = \JSE(T)^*$.
\item Given an increasing tableau $T$ of partition shape define $T^\flat = \JNW(T^*)$.
\end{enumerate}
\end{defn}

\begin{prop} \label{P:decreasing to increasing} The map $T\mapsto T^\sharp$ is a bijection
$\Dec_\la\to\Inc_\la$ with inverse $T\mapsto T^\flat$. Moreover
$T^\sharp=\JNW(T^*)$ and $T^\flat=\JSE(T)^*$.
\end{prop}
\begin{proof} 
It follows from Proposition \ref{P:Kjdt Pieri} that infusion sends a CSS tableau to a CR tableau whose shape is the 180-degree rotation, and vice versa. This implies that $T\mapsto T^\sharp$ and $T\mapsto T^\flat$ are shape-preserving. They are mutually inverse since $*$ is involutive and infusion is involutive.
The alternate descriptions of the maps hold since switching commutes with 180 degree rotation.
\end{proof}

\begin{ex} \label{X:Krect} Let $T$ be the following decreasing tableau $T$.
We illustrate how to compute $T^\sharp$.
The first step is to anti-rectify $T$ with 
respect to
the rectification order $U$. 
$\JSE(T)$ is the result of 
the anti-rectification.
\[
T = 
\begin{ytableau}
9 & 8 & 7 & 5 & 3\\
7 & 5 & 4 & 3\\
4 & 2 & 1 \\
3\\
1\\
\end{ytableau}
\quad 
U = \ytableaushort{\none\none\none\none,\none\none\none1,\none\none25,368a,479b}\quad
\JSE(T) = 
\begin{ytableau}
\none & \none & \none & \none & 7\\
\none & \none & \none & \none & 5\\
\none & \none & 8 & 7 & 3 \\
\none & 9 & 5 & 4 & 2\\
9 & 7 & 3 & 2 & 1\\
\end{ytableau}
\]

Then we rotate $\JSE(T)$ and obtain $T^\sharp$:
$$
T^\sharp = 
\begin{ytableau}
1 & 2 & 3 & 7 & 9\\
2 & 4 & 5 & 9\\
3 & 7 & 8 \\
5\\
7\\
\end{ytableau}
$$

Next, we illustrate how to compute $(T^\sharp)^\flat$.
We rotate $T^\sharp$ and get $(T^\sharp)^*$.
Then we rectify $(T^\sharp)^*$ using $U'$.
The result $\JNW((T^\sharp)^*)$ is $(T^\sharp)^\flat$.
\[
(T^\sharp)^* = 
\begin{ytableau}
\none & \none & \none & \none & 7\\
\none & \none & \none & \none & 5\\
\none & \none & 8 & 7 & 3 \\
\none & 9 & 5 & 4 & 2\\
9 & 7 & 3 & 2 & 1\\
\end{ytableau}
\quad 
U' = 
\begin{ytableau}
1 & 5 & 8 & a\\
2 & 6 & 9 & b\\
3 & 7\\
4\\
\end{ytableau}
\quad 
\JNW((T^\sharp)^*) = 
\begin{ytableau}
9 & 8 & 7 & 5 & 3\\
7 & 5 & 4 & 3\\
4 & 2 & 1 \\
3\\
1\\
\end{ytableau}
\]
Readers can check $(T^\sharp)^\flat = T$ 
in this case. 
\end{ex}

\begin{rem} 
\begin{enumerate}
\item Unlike jeu-de-taquin for semistandard tableaux, $\JNW_S(T)$ may depend on the rectification order $S$ \cite[Ex. 1.3]{TY}. Thus for the well-definedness of $\JNW(T)$, it is necessary to specify $S$.
\item The $\la$ and $\mu$ defining a skew shape $\la/\mu$ are not unique; one may add several rows and columns to the top and left of the diagrams of $\mu$ and $\la$ simultaneously and get the same difference of partition diagrams. Using Proposition \ref{P:Kjdt Pieri} it can be shown that $\JNW(T)$ depends only on the set of boxes in $\la/\mu$ and not on the pair $(\la,\mu)$.
\item This definition of $\JSE(T)$ was used in \cite{ReY}.
\end{enumerate}
\end{rem}

\newcommand{\da}{
\begin{tikzpicture}[scale=.5]
\draw[black,thick] (0,0)--(8,0)--(8,5)--(0,5)--(0,0);
\draw[black,thick] (1,0)--(1,1)--(3,1)--(3,3)--(6,3)--(6,4)--(8,4);
\node (T) at (2,4) {$T^\sharp$ incr.};
\node (CSS) at (6,2) {CSS};
\end{tikzpicture}
}

\newcommand{\db}{
\begin{tikzpicture}[scale=.5]
\draw[black,thick] (0,0)--(8,0)--(8,5)--(0,5)--(0,0);
\draw[black,thick] (0,1)--(2,1)--(2,2)--(5,2)--(5,4)--(7,4)--(7,5)--(8,5);
\node (TJSE) at (6,1) {$T^*$ incr.};
\node (CSS) at (2,4) {CR};
\end{tikzpicture}
}

\newcommand{\dc}{
\begin{tikzpicture}[scale=.5]
\draw[black,thick] (0,0)--(8,0)--(8,5)--(0,5)--(0,0);
\draw[black,thick] (0,1)--(2,1)--(2,2)--(5,2)--(5,4)--(7,4)--(7,5)--(8,5);
\node (Tstar) at (6,1) {$\JSE(T)$ decr.};
\node (CSS) at (2,4) {CSS};
\end{tikzpicture}
}

\newcommand{\dd}{
\begin{tikzpicture}[scale=.5]
\draw[black,thick] (0,0)--(8,0)--(8,5)--(0,5)--(0,0);
\draw[black,thick] (1,0)--(1,1)--(3,1)--(3,3)--(6,3)--(6,4)--(8,4);
\node (Tflat) at (2,4) {$T$ decr.};
\node (CSS) at (6,2) {CR};
\end{tikzpicture}
}

\newsavebox{\myboxa}
\newsavebox{\myboxb}
\newsavebox{\myboxc}
\newsavebox{\myboxd}

\tikzstyle{arrow} = [thick,<->,>=stealth]

\newcommand{\twobytwo}[4]{ %
\sbox{\myboxa}{#1}
\sbox{\myboxb}{#2}
\sbox{\myboxc}{#3}
\sbox{\myboxd}{#4}
\begin{tikzpicture}
\node (a) {\usebox{\myboxd}};
\node (b) [right = 1.5cm of a] {\usebox{\myboxc}};
\node (c) [below = 1.5cm of a] {\usebox{\myboxb}};
\node (d) [right = 1.5cm of c] {\usebox{\myboxa}};
\draw [arrow] (a) -- (b) node[midway,yshift=10pt] {\text{infusion}}; 
\draw [arrow] (a) -- (c) node[midway,xshift=-10pt] {*}; 
\draw [arrow] (b) -- (d) node[midway,xshift=10pt] {*};
\draw [arrow] (c) -- (d) node[midway,yshift=-10pt] {\text{infusion}}; 

\end{tikzpicture}}

\begin{figure}
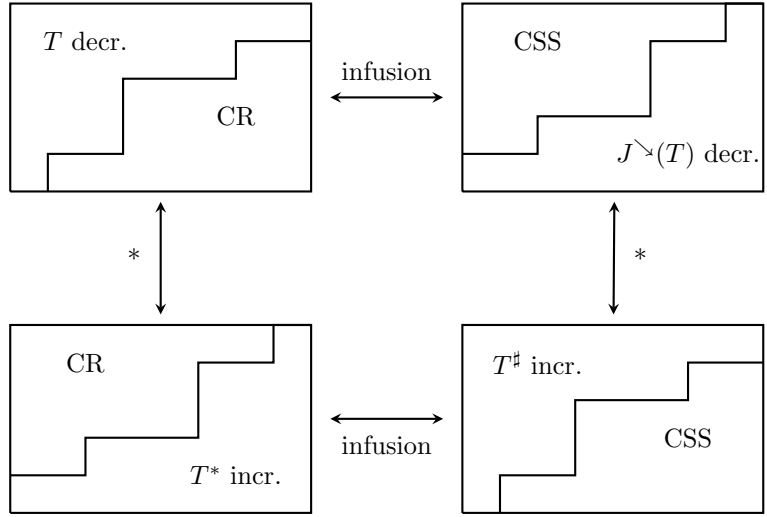

\twobytwo{\da}{\db}{\dc}{\dd}
\caption{From decreasing to increasing: $T\mapsto T^\sharp$}
\end{figure}

\subsubsection{Keys defined via Kjdt for increasing/decreasing tableaux}
\label{SS:Kjdt keys}
For a tableau $T$ let $T_{\le j}$ be the tableau consisting of the first $j$ columns of $T$. $T_{\ge j}$ is defined similarly.

Let $T$ be an increasing (resp. decreasing) tableau of partition shape.
The \emph{left key} $K_-(T)$ of the increasing tableau $T$ is defined by the condition that its $j$-th column is equal to the first column of $\JSE(T_{\le j})$ for all $j$.
The \emph{right key} $K_+(T)$ of the decreasing tableau $T$ is defined by the condition that its $j$-th column is equal to the last column of $\JSE(T_{\ge j})$ for all $j$.
We note that $\JSE$ is defined using a different reverse rectification order for increasing versus decreasing tableaux. However, it is ultimately shown that these keys are independent of the rectification order; 
see Proposition \ref{P:Kjdt right key} for decreasing tableaux.

\subsection{Hecke insertion}
\label{SS:Hecke insertion}
Let $T$ be a decreasing tableau of partition shape $\la$ and $x$ 
a positive integer. The (column) \emph{Hecke insertion} of $x$ into $T$ is defined as follows. It maps the pair $(x,T)$ to a triple 
$(P,c, \alpha)$.
$P$ will be a decreasing tableau which either has shape $\la$ in which case we set $\alpha=0$ or differs by adding a single box, in which case we set $\alpha=1$. We write $P = (x\HD T)$.
$c$ will be the box of $P$ where the algorithm ends.

The algorithm first inserts $x$ into column 1 of $T$.
This may output a number. 
If so, the output number is then inserted to the next column.
The algorithm repeats until an insertion
to a column has no output.
To describe the insertion of $x$ to a column of $T$,
we consider two cases:
\begin{enumerate}
\item[Case 1:] $x$ is less than or equal 
to all entries in this column.
Then the algorithm makes no output. 
In addition,
it appends $x$ to the bottom of this column as long 
as the result is a decreasing tableau.
$c$ is set to be this newly appended box.
Otherwise the column is unchanged and 
$c$ is set to be the rightmost box in the row
that contains the bottom entry of this column.
\item[Case 2:] Otherwise let $y$ be the smallest value in
this column such that $y < x$.
Then the algorithm outputs $y$ from the column.
In addition it replaces $y$ by $x$ as long 
as the result is a decreasing tableau.
Otherwise the column remains unchanged.
\end{enumerate}
We use the word ``contraction" when $\alpha=0$.

This algorithm has an inverse called
reverse (column) Hecke insertion,
which maps a triple $(P,c, \alpha)$ 
to a pair $(T,x)$.
Here, $P$ is a decreasing tableau of partition
shape. 
$c$ is an entry on $P$ that is at the end of 
its row and its column. 
$\alpha$ is 0 or 1.

The algorithm behaves in the following way:
First let $y$ be the number in the box $c$.
If $\alpha = 1$, we remove this box.
If $\alpha = 0$, we do not remove it.
In either case, the algorithm ``reverse inserts'' 
$y$ into the previous column.
When a value $y$ is ``reverse inserted''
into a column, the algorithm finds the largest
$y'$ in that column such that $y' > y$. 
It replaces $y'$ by $y$ as long 
as the result is a decreasing tableau.
Otherwise it does nothing to the column.
In either case, $y'$ is reverse inserted
into the previous column. 
If there is no column on the left, 
the algorithm lets $x = y'$ and $T$ is 
the resulting tableau. Then it terminates. 

\begin{lemma} \cite{BKSTY}
\label{Hecke insertion are inverses}
Hecke insertion and reverse Hecke insertion
are inverses of each other.
\end{lemma}

\subsubsection{Pieri property of Hecke insertion}

Hecke insertion has the following Pieri property:
\begin{lemma}\cite[Lemma 2]{BKSTY}
\label{Pieri property of Hecke insertion}
Let $T$ be a decreasing tableau.
Let $x_1, x_2$ be two positive integers. 
Hecke insert $x_1$ into $T$ with result $(T_1,c_1,\alpha_1)$ 
and then Hecke insert $x_2$ into $T_1$, with result $(T_2,c_2,\alpha_2)$.
Then $c_2$ is strictly to the right of $c_1$ if and only if
$x_1 < x_2$.
\end{lemma}

\subsubsection{Set-valued recording tableaux for Hecke insertion and Hecke RSK}
\label{SS:SV recording}
Let $\T := \bigsqcup_\la (\Dec_\la \times \RSVT_\la)$.
Now we describe the bijections between
$\C$ and $\T$.
First, we recursively define 
$\Insert: \C \mapsto \T$.
\begin{defn}
Take $(a,i) \in \C$.
If $a$ is the empty string,
then $\Insert(a,i)$ is the pair
of two empty tableaux.
Now assume $a$ has positive length.
Let $a = a'x$, $i = i'y$,
where $x, y$ are positive numbers.
Now let $(P', Q') = \Insert(a', i')$.
We Hecke insert $x$ into $P'$ and 
get $(P,c,\alpha)$.
If $\alpha = 1$,
we append $y$ to $Q'$ at the corresponding
position of $c$.
Otherwise, 
we add $y$ to the entry in $Q'$ that
corresponds to $c$.
We let $Q$ be the resulting RSVT.
Then $\Insert(a,i) := (P,Q)$.
\end{defn}

$\Insert$ is well-defined by Lemma \ref{Pieri property of Hecke insertion}.

Now we recursively define the map
$\RevInsert: \T \mapsto \C$.
\begin{defn}
Take $(P,Q) \in \T$.
If $P$ is the empty tableau,
$\RevInsert(P,Q)$ is the pair of two
empty words.
Now assume $P$ is non-empty. 
Let $y$ be the smallest number in $Q$.
We pick the rightmost $y$, 
and remove this number from $Q$.
Let $Q'$ be the resulting RSVT. 
If this $y$ is the only number in its entry,
we set $\alpha = 1$. 
Otherwise, $\alpha = 0$.
Then we invoke reverse Hecke insertion on the 
corresponding entry in $P$ with $\alpha$.
Let $(P', x)$ be the output.
Finally, we let $(a', i') = \RevInsert(P', Q')$.
Then $\RevInsert(P, Q) := (ax, iy)$.
\end{defn}

$\RevInsert$ is similarly well-defined.

\begin{lemma}
$\Insert$ and $\RevInsert$ are inverses.
\end{lemma}
\begin{proof}
Follows from Lemmas \ref{Hecke insertion are inverses} and 
\ref{Pieri property of Hecke insertion}.
\end{proof}

Also, these two maps have the following 
property:
\begin{prop} \label{P:Hecke insertion Kequiv and weight}
Take $(a,i) \in \C$ and let $(P,Q) = \Insert(a,i)$.
\begin{align}
\label{E:insert respects permutations}
\rev(a) &\Kequiv P \\
\label{E:insert is weight preserving}
\wt(Q)&=\wt(i)
\end{align}
\end{prop}
\begin{proof} Equation \eqref{E:insert is weight preserving} holds by definition.
The relation \eqref{E:insert respects permutations} holds by \cite[Thm. 6.2]{BS}.
\end{proof}

\section{RSVT rule 
for Lascoux Polynomials}
\label{S:formula for Lascoux polynomials}
In this section, 
we give a combinatorial rule for Lascoux polynomials
involving tableaux.
Let $T$ be a RSVT.
Let $L(T)$ be the RSSYT obtained
by picking the largest number in each entry.
Then we have
\begin{thm} For any composition $\alpha$
\label{set-valued_reverse_SSYT_rule}
\begin{align}
\fLb_\alpha
= \sum_{K_-(L(T)) \leq \alpha} 
\beta^{\extra(T)} x^{\wt(T)}
\end{align}
where $T$ runs over the RSVT of shape $\alpha^+$.
\end{thm}
\begin{ex}
The following RSVTs 
contribute to $\fLb_{(1,0,2)}$:\\
\begin{ytableau}
2 & 1\\
1 
\end{ytableau}
\\[1mm]
\begin{ytableau}
2 & 2\\
1 
\end{ytableau}
\quad
\begin{ytableau}
2 & 21\\
1 
\end{ytableau}
\\[1mm]
\begin{ytableau}
3 & 1\\
1 
\end{ytableau}
\quad
\begin{ytableau}
32 & 1\\
1 
\end{ytableau}
\\[1mm]
\begin{ytableau}
3 & 2\\
1 
\end{ytableau}
\quad
\begin{ytableau}
3 & 21\\
1 
\end{ytableau}
\quad
\begin{ytableau}
32 & 2\\
1 
\end{ytableau}
\quad
\begin{ytableau}
32 & 21\\
1 
\end{ytableau}
\\[1mm]
\begin{ytableau}
3 & 3\\
1 
\end{ytableau}
\quad
\begin{ytableau}
3 & 31\\
1 
\end{ytableau}
\quad
\begin{ytableau}
3 & 32\\
1 
\end{ytableau}
\quad
\begin{ytableau}
3 & 321\\
1 
\end{ytableau}
\vskip2mm

Thus, we may write $\fLb_{(1,0,2)}$ as
\begin{equation*}
\begin{split}
&x_1^2x_2 + x_1x_2^2 + x_1^2x_3 
+ x_1x_2x_3 + x_1x_3^2 \\
+ & \beta(x_1^2x_2^2 + 2x_1^2x_2x_3 + 
x_1x_2^2x_3 + x_1^2x_3^2 + x_1x_2x_3^2)\\
+ & \beta^2(x_1^2x_2^2x_3 + x_1^2x_2x_3^2)
\end{split}
\end{equation*}
\end{ex}

Next, we want to rewrite
Theorem \ref{set-valued_reverse_SSYT_rule}
as a rule involving RSSYT, instead of RSVT. 
We need a definition:
\begin{defn}
Fix a RSSYT $T$ of shape $\lambda$. Define $WT(T)$ by
$$
x^{WT(T)} =
\sum_{\substack{T' \in \RSVT_\la \\ L(T') = T}} 
\beta^{\extra(T')} x^{\wt(T')}.
$$
\end{defn}
Then Theorem \ref{set-valued_reverse_SSYT_rule}
can be rewritten as:
\begin{thm} For any composition $\alpha$
\label{reverse_SSYT_rule}
$$
\fLb_\alpha
= \sum_{\substack{T\in\RSVT_{\alpha^+} \\ K_-(T) \leq \alpha}} x^{WT(T)}
$$
\end{thm}
It is clear that Theorem $\ref{reverse_SSYT_rule}$
and Theorem $\ref{set-valued_reverse_SSYT_rule}$
are equivalent.
Readers may still insist that Theorem $\ref{reverse_SSYT_rule}$
involves RSVTs,
since they appear in how we defined $x^{WT(T)}$.
Well, the following lemma resolves this issue:
\begin{lem} \label{L:WT} For any RSSYT of shape $\la$
\begin{align}
x^{WT(T)} = x^{\wt(T)}\prod_{(s,k)}(1 +\beta x_k)
\end{align}
where $(s,k)$ runs over pairs such that $s$ is a box in $\la$,
$k$ is less than the value of $T$ in that box, and replacing the $s$-th entry of $T$ by $k$ results in a RSSYT.
\end{lem}

\begin{proof}[Proof of Lemma \ref{L:WT}]
Consider the following way of turning
$T$ into a RSVT in $L^{-1}(T)$.
Let $a$ be an entry in $T$.
Let $b$ be the entry on its right and 
$b = 1$ if such an entry does not exist.
Let $c$ be the entry below $a$ and 
$c = 0$ if such an entry does not exist.
We turn $a$ into $\{ a \}$, 
and then add some numbers to this set. 
We may add any $k$ such that $a > k$, $k > c$
and $k \geq b$. 
Not adding this $k$ will contribute $1$
and adding this $k$ will contribute an $\beta x_k$.
Thus, each such $k$ contributes $(1 +\beta x_k)$.
Clearly, the choices are independent 
and any element in $L^{-1}(T)$ can be obtained
this way. 
\end{proof}

In the rest of this section,
we show Theorem
$\ref{reverse_SSYT_rule}$.
We will only use RSSYTs.
We only need to prove the sum in 
Theorem $\ref{reverse_SSYT_rule}$
satisfies the recursion of Lascoux polynomials.
Now, we fix an $i$ throughout the rest of this section.

\subsection{Partitioning RSSYTs}

Let $T$ be a RSSYT.
We classify its $i$ and $i+1$ 
into 3 categories: ``ignorable",
``frozen", and ``free".
First, we find all pairs of $i+1$ 
and $i$ that appear in 
the same column.
We pair them 
and say they are ``ignorable".
Next, we find non-ignorable $i$ and $i+1$ 
such that:
\begin{enumerate}
\item $i$ is on the left of $i+1$.
\item Any column between them
must have an ignorable pairs.
\end{enumerate}
We pair them and say they are ``frozen".
Other non-ignorable $i$ and $i+1$ are
called ``free".
\begin{ex}
When $i = 3$, 
consider the following RSSYT:
\begin{align*}
\begin{ytableau}
6 & 6 & 6 & 6 & \textcolor{blue}{4} & 4\\
5 & \textcolor{red}{4} & 3 & \textcolor{blue}{3}\\
4 & \textcolor{red}{3}
\end{ytableau} 
\end{align*}

The red entries are ignorable
and blue entries are frozen.
Other 3 and 4 are free.
\end{ex}

Based on this labelling, 
we may partition RSSYTs into families.
\begin{defn}
A \emph{family} is an equivalence class 
under the transitive closure of the 
following: two RSSYTs are related if they differ by 
changing a single $i$ into an $i+1$ (or vice versa) 
where the changed letters are free in both tableaux.
\end{defn}

\begin{ex}
Consider the reverse tableau 
in the previous example. 
Its family also includes:
\begin{align*}
\begin{ytableau}
6 & 6 & 6 & 6 & \textcolor{blue}{4} & 3\\
5 & \textcolor{red}{4} & 3 & \textcolor{blue}{3}\\
3 & \textcolor{red}{3}
\end{ytableau}
\quad
\begin{ytableau}
6 & 6 & 6 & 6 & \textcolor{blue}{4} & 3\\
5 & \textcolor{red}{4} & 3 & \textcolor{blue}{3}\\
4 & \textcolor{red}{3}
\end{ytableau}
\quad
\begin{ytableau}
6 & 6 & 6 & 6 & \textcolor{blue}{4} & 3\\
5 & \textcolor{red}{4} & 4 & \textcolor{blue}{3}\\
4 & \textcolor{red}{3}
\end{ytableau}
\end{align*}

\begin{align*}
\begin{ytableau}
6 & 6 & 6 & 6 & \textcolor{blue}{4} & 4\\
5 & \textcolor{red}{4} & 3 & \textcolor{blue}{3}\\
3 & \textcolor{red}{3}
\end{ytableau}
\quad
\begin{ytableau}
6 & 6 & 6 & 6 & \textcolor{blue}{4} & 4\\
5 & \textcolor{red}{4} & 4 & \textcolor{blue}{3}\\
4 & \textcolor{red}{3}
\end{ytableau}
\end{align*}

However, the following is in another family:
\begin{align*}
\begin{ytableau}
6 & 6 & 6 & 6 & \textcolor{blue}{4} & 3\\
5 & \textcolor{red}{4} & \textcolor{blue}{4} & \textcolor{blue}{3}\\
\textcolor{blue}{3} & \textcolor{red}{3}
\end{ytableau}    
\end{align*}
\end{ex}

Given an RSSYT, 
how can we enumerate its family?
Clearly, we can only change its free entries.
We also need to make sure they are still
free after our changes.
In other words, 
assume $a$ and $b$ are two
free entries.
If $a$ is on the left of $b$ 
and all columns between
them have ignorable pairs,
then we cannot change $a$ into $i$
and $b$ into $i + 1$.
This criterion leads 
to the following definition.
\begin{defn}
Let $T$ be a RSSYT.
We partition its free $i$ 
and $i+1$ into ``blocks".
Two entries are in the same block iff
all columns between them 
have ignorable pairs.
\end{defn}
Thus, to enumerate the family 
of a RSSYT $T$,
we just replace entries
in each block by a weakly decreasing 
(from left to right)
sequence of $i$ and $i+1$.
The reader may check the enumeration of 
the family in the previous example.

\subsection{Families and left keys}
\label{SS:families and left keys}
This subsection aims to describe 
the left keys of a family.
This idea is formalized 
in the following lemma:

\begin{lem}
\label{family_and_left_key}
Let $\mathcal{F}$ be a family.
Then its elements can have at most
2 different left keys. 
If they all have the same left key 
$\gamma$,
then $\gamma_i \geq \gamma_{i+1}$.

If they have two different left keys,
then they must be $\gamma$ 
and $s_i \gamma$, 
where $\gamma_i > \gamma_{i+1}$.
In this case, we also have:
\begin{enumerate}
\item $T \in \mathcal{F}$ has left key 
$\gamma$ iff
$T$'s leftmost block only has $i$.
\item All columns before
the first block must have ignorable pairs.
\end{enumerate}
\end{lem}

Before proving the lemma, 
we need to introduce an algorithm
that computes the left key.
The algorithm is introduced in section 5
of \cite{Willis}.
Here we describe this algorithm 
in a slightly different way. 
\begin{defn}
Given two columns $C_1, C_2$ such that 
$C_1C_2$ is a RSSYT, we define the 
column $C_1 \triangleleft C_2$ as follows.
Assume $C_2 = \{ a_1 < a_2 < \dots < a_m\}$.
We find the smallest $b_1$ in $C_1$
such that $b_1 \geq a_1$.
Then we find the smallest $b_2$ in $C_1$
such that $b_2 \geq a_2$ and $b_2 > b_1$.
Similarly, we find $b_3, \dots b_m$.
Let $C_1 \triangleleft C_2=\{b_1<b_2<\dotsm<b_m\}$.

More generally suppose $C_1,C_2,\dotsc,C_k$ 
are the columns of a RSSYT. 
Observe that the following expression is well-defined
when $j \le k$
$$
C_j \triangleleft \dotsm \triangleleft C_k := 
C_j \triangleleft (C_{j+1} \triangleleft \dotsm \triangleleft C_k)
$$
where the base case is
$$
C_k \triangleleft \dotsm \triangleleft C_k := C_k
$$
\end{defn}
With this definition,
we may compute column $k$ of $K_-(T)$,
where $T$ is a RSSYT. 
Let the first $k$ columns of $T$
be $C_1, \dots, C_k$.
Then column $k$ of $K_-(T)$ 
is by definition $C_1 \triangleleft \dots \triangleleft C_k$.

To study this algorithm,
we need to classify columns of $T$.
Each column can be labeled as follows:
\begin{enumerate}
\item Type 1 column: It has neither $i$
nor $i+1$.
\item Type 2 column: It has $i$ but no
$i+1$.
\item Type 3 column: It has $i+1$ but no
$i$.
\item Type 4 column: It has both $i$ and 
$i+1$.
\end{enumerate}

Now we make several observations.
\begin{lemma}
\label{left_key_algorithm_property_1}
If $C_1$ has type 4 and $C_2$ does not
have type 3,
then $C_1 \triangleleft C_2$ cannot have type 3.
\end{lemma}
\begin{proof}
Assume $C_1 \triangleleft C_2$ has type 3.
Then we must pick $i+1$ in $C_1$ for
some $m$ in $C_2$.
Moreover, $i$ in $C_1$ is never picked.
Thus, $m$ must be $i+1$ and $C_2$ cannot
have $i$.
$C_2$ has type 3, contradiction.
\end{proof}

\begin{lemma}
\label{left_key_algorithm_property_2}
Let $T$ be a RSSYT
with no free $i+1$.
Assume $\gamma = K_-(T)$.
Then $\gamma_i \geq \gamma_{i+1}$.
\end{lemma}
\begin{proof}
Let $C_1, C_2, \dots$ be columns of $T$.
Consider column $k$ of $K_-(T)$.
We only need to prove it cannot have
type 3.

Suppose $C_1, \dots, C_k$ all have type 4.
Then Lemma \ref{left_key_algorithm_property_1}
guarantees $C_1 \triangleleft \dotsm \triangleleft C_{k}$ 
cannot have type 3. 
Otherwise, we can find $j \leq k$ such that
$C_1, \dots, C_{j-1}$ have type 4
and $C_j$ does not have type 4.
Since $T$ has no free $i+1$,
$C_j$ must have type 1 or 2.
Then $C_j \triangleleft \dotsm \triangleleft C_{k}$ 
must also have type 1 or 2.
By Lemma \ref{left_key_algorithm_property_1},
$C_1, \dots, C_{j-1}$ cannot turn it 
into type 3. 
\end{proof}

\begin{lemma}
\label{left_key_algorithm_property_3}
Assume $C_2$ has type 2.
We change its $i$ into $i+1$
and obtain $C_2'$.
Assume $C_1C_2'$ 
is a RSSYT.
Then, 
\begin{enumerate}
\item If $C_1$ has type 4,
then $C_1 \triangleleft C_2 = 
C_1 \triangleleft C_2'$,
or $C_1 \triangleleft C_2'$ is obtained from
$C_1 \triangleleft C_2$ by changing an $i$
into $i+1$.
\item If $C_1$ has type 1 or 3,
then $C_1 \triangleleft C_2 = 
C_1 \triangleleft C_2'$.
\end{enumerate}
\end{lemma}
\begin{proof}
We do a case study based on the type 
of $C_1$.
\begin{enumerate}
\item Assume $C_1$ has type 4.
When we consider $i$ in $C_2$,
there are 3 possibilities:
The $i$ in $C_1$ is picked; 
or a number larger than
it is picked; 
or the $i$ is still available. 

In the first 2 cases, 
clearly this $i$ in $C_2$ behaves
as if it is an $i+1$.
Then $C_1 \triangleleft C_2 = 
C_1 \triangleleft C_2'$.
In the last case,
$i$ in $C_2$ picks $i$,
and $i+1$ in $C_2'$ picks $i+1$.
Our claim is clear. 
\item Assume $C_1$ has type 1 or 3.
Clearly the $i$ in $C_2$ behaves
as if it is an $i+1$,
so $C_1 \triangleleft C_2 = 
C_1 \triangleleft C_2'$.
\end{enumerate}
\end{proof}

\begin{lemma}
\label{left_key_algorithm_property_4}
Let $T$ be a RSSYT.
Assume column $j$ of $T$ has a free $i$,
which is the leftmost free $i$ in its block.
We change this $i$ into $i+1$ and get $T'$.
If $\gamma = K_-(T)$, 
then $K_-(T') = \gamma$ or $s_i \gamma$.
Moreover, if the latter case happens,
we must have: 
\begin{enumerate}
\item The $i$ we changed is in the 
leftmost block of $T$.
\item Each of column $1, \dots, j-1$ of $T$
has ignorable pairs. 
\end{enumerate}
\end{lemma}
\begin{proof}
Let $C_1, C_2, \dots$ be the columns of $T$.
Let $D_1, D_2, \dots$ be the columns of $T'$.
Consider column $k$ of $K_-(T)$ 
and $K_-(T')$.
If $k < j$,
then clearly they are the same.
Now assume $k > j$.
Let $C = C_{j+1} \triangleleft \dotsm \triangleleft C_{k}$.
Because the $i$ in column $j$ is free,
we know that $C_{j+1}, \dots, C_k$ all have type 4,
or the leftmost non-type-4 column among
them has type 1 or 2.
Similar to the proof of Lemma
\ref{left_key_algorithm_property_2},
$C$ cannot have type 3.
Next, we compare $C_j \triangleleft C$
and $D_j \triangleleft C$.
If $i$ in $C_j$ is picked by $x$ in $C$,
then this $x$ will pick $i+1$ in $D_j$.
Thus, $D_j \triangleleft C$ is obtained
by changing $i$ in $C_j \triangleleft C$
into $i+1$.
If $i$ in $C_j$ is not picked,
the $i+1$ in $D_j$ will not be picked.
Then $C_j \triangleleft C = D_j \triangleleft C$.

Consequently, if $k \geq j$,
$C_j \triangleleft \dotsm \triangleleft C_k$
agrees with 
$D_j \triangleleft \dotsm \triangleleft D_k$,
or the latter differs from the former
by changing an $i$ into $i+1$.
In Lemma \ref{left_key_algorithm_property_3},
we showed this difference might be preserved
or corrected by type 4 columns.
If $C_1, \dots, C_{j-1}$ all have type 4,
then we know column $k$ of $K_-(T)$ agrees
with column $k$ of $K_-(T')$, 
or the latter differs from the former
by changing an $i$ into $i+1$.
Otherwise, 
we let $l$ be the largest such that 
$l < j$ and $C_l$ does not have type 4.
Since the $i$ in column $j$ of $T$ is the 
leftmost $i$ in its block,
$C_l$ must have type 1 or 3.
By Lemma \ref{left_key_algorithm_property_3},
$$
C_l \triangleleft \dotsm \triangleleft C_k 
= D_l \triangleleft \dotsm \triangleleft D_k
$$

Thus, each column of $K_-(T')$ either
agrees with the corresponding 
column in $K_-(T)$,
or differs by changing an $i$ into $i+1$.
Since $K_-(T')$ is a key,
we have $K_-(T') = \gamma$ or $s_i \gamma$.
In the latter case
we know $C_1, \dots, C_{j-1}$
have type 4.
Our claims are immediate.

\end{proof}

Now we may prove Lemma \ref{family_and_left_key}.

\begin{proof}
First pick $T$ from 
$\mathcal{F}$ that has no free $i+1$.
Assume $\gamma = K_-(T)$.
By Lemma \ref{left_key_algorithm_property_2},
$\gamma_i \geq \gamma_{i+1}$.

Then we enumerate other elements 
in $\mathcal{F}$ by changing
free $i$ in $T$ into $i+1$.
As long as we do not change the first block,
the left key will still be $\gamma$.
Once we change the first $i$ in the 
first block,
the left key might be fixed,
or turned into $s_i \gamma$.
The latter case is possible only when
all columns before the first blocks
have ignorable pairs.
After that, no matter which $i$ we change, 
the left key will be fixed. 
\end{proof} 
\subsection{$\pi_i$ and $\pi_i^K$}
In this subsection, 
we derive some basic facts 
about $\pi_i$ and $\pi_i^K$.
Define $X_i = x_i(1 + \beta x_{i+1})$ 
and $X_{i+1} = x_{i+1}(1 + \beta x_i)$.
Then we have
\begin{enumerate}
\item $s_i (X_i) = X_{i+1}$
\item $\pi_i(f) = \partial_i(x_if)$ and 
$\pi_i^K(f) = \partial_i(X_if)$
\item $\partial_i(X_i) = \partial_i(x_i) = 1$.
\end{enumerate}

The following lemma describes how $\partial_i$
acts on a product of several $x_i$ and $X_i$:
\begin{lemma}
Assume we have $u_1, \dots, u_n$,
where each $u_j$ is either $x_i$ or $X_i$.
Then 
$$
\partial_i(u_1 \dots u_n) = 
\sum_{j = 1}^n s_i(u_1 \dots u_{j-1})u_{j+1}\dots u_n
$$
\end{lemma}

For instance, 
$$
\partial_i(x_iX_ix_iX_i) = 
X_ix_iX_i + x_{i+1}x_iX_i + 
x_{i+1}X_{i+1}X_i + x_{i+1}X_{i+1}x_{i+1}
$$

\begin{proof}
Notice:
\begin{equation*}
\begin{split}
\partial_i(u_1 \dots u_n) & = 
\partial_i(u_1) u_2 \dots u_n 
+ s_i(u_1)\partial_i(u_2 \dots u_n) \\
& = u_2 \dots u_n 
+ s_i(u_1)\partial_i(u_2 \dots u_n)   
\end{split}
\end{equation*}

Then the proof is finished by induction.
\end{proof}

\begin{cor}
\label{pi_comb}
Assume we have $u_1, \dots, u_n$,
where each $u_j$ is either $x_i$ or $X_i$.
Then 
\begin{align}
\pi_i(u_1 \dots u_n) = 
u_1 \dots u_n
+ x_{i+1}\sum_{j = 1}^n 
s_i(u_1 \dots u_{j-1})u_{j+1}\dots u_n
\end{align}
\begin{align}
\pi_i^K(u_1 \dots u_n) = 
u_1 \dots u_n
+ X_{i+1}\sum_{j = 1}^n 
s_i(u_1 \dots u_{j-1})u_{j+1}\dots u_n
\end{align}
\end{cor}

\subsection{$x^{WT(T)}$ and Family}
In this subsection, 
we investigate how $x^{WT(T)}$ works 
and how it 
changes within a family.
More explicitly, 
the goal is to understand:
$\sum_{T \in \mathcal{F}} x^{WT(T)}$
where $\mathcal{F}$ is a family.

The first step
is to understand what governs the
power of $(1 + \beta x_j)$ in $x^{WT(T)}$.
Based on our definition, 
each row can have 
at most one entry 
that contributes $(1 + \beta x_j)$ for a fixed j.
How is it determined whether a row has 
such a contributor?
The following lemma answers this question.
To make it concise,
we adopt the following convention throughout
the rest of this section:
a 0 is appended below each column 
in a RSSYT.

\begin{lemma}
A row has an entry that contributes $(1 
+ \beta x_j)$
iff we can find an entry $j'$ on this row 
such that:
\begin{enumerate}
\item $j' > j$
\item The entry below $j'$ is less than $j$. 
\end{enumerate}
\end{lemma}
\begin{proof}
Assume an entry $m$ contributes $1 + \beta x_j$.
Then clearly $m > j$ and 
the entry below $m$ is less than $j$.
The row of $m$ clearly satisfies
the requirement.

Conversely, assume a row has $j'$
that satisfies the two requirements.
Moreover, we pick the rightmost $j'$ 
among all such $j'$ on this row.
Then the entry to the right of $j'$
either does not exist or is at most $j$.
Changing this $j'$ to $j$ will make $T$
a valid anti-SSYT.
Thus, this entry contributes 
$(1 + \beta x_j)$.
\end{proof}

With this lemma,
we may ascribe contributions of 
$(1 + \beta x_j)$ to rows, 
instead of entries.
However, we would like to ascribe
contributions of 
$(1  + \beta x_i)$ and $(1 + 
\beta x_{i+1})$
to specific entries, 
but the rule is different from 
our previous criterion. 
If a row contributes $(1 + \beta x_i)$, 
then we may find the leftmost entry 
on this row satisfying:
\begin{enumerate}
\item It is larger than $i$.
\item The entry below it is less than $i$.
\end{enumerate}
We say this entry 
contributes an $(1 + \beta x_i)$.
Similarly, 
if a row contributes 
$(1 + \beta x_{i+1})$, 
then we may find the rightmost entry 
on the row below satisfying:
\begin{enumerate}
\item It is less than $i+1$.
\item The entry above it is larger than $i+1$.
\end{enumerate}
We say this entry contributes an $(1 + \beta x_{i+1})$.
To illustrate our new ``contribution system",
consider the following example:

\begin{ex}

$$
\begin{ytableau}
6 & 6 & 6 & 6 & \textcolor{blue}{4} & 4\\
5 & 4 & \textcolor{blue}{4} & \textcolor{red}{3} & 0 & 0\\
\textcolor{blue}{4} & 3 & 0 & 0\\
0 & 0
\end{ytableau}
$$
When $i = 3$,
each blue 4 contributes 
$x_4(1 + \beta x_3)$.
The red 3 contributes 
$x_3(1 + \beta x_4)$.
\end{ex}

Now we fix an arbitrary family $\mathcal{F}$
throughout this subsection.
Take any $T \in \mathcal{F}$.
Let $m$ be the number of blocks in $T$.
Then we may break $x^{WT(T)}$ into a product:
$$
x^{WT(T)} = g^Tf_1^T\dots f_m^T
$$
Here, 
$f_j^T$ is the contribution of the 
$j^{th}$ block in $T$ 
from left to right.
$g^T$ contains the contribution 
of $x_i$, $x_{i+1}$, $(1 + \beta x_i)$ 
and $(1 + \beta x_{i+1})$ 
from all other entries.
It also contains powers of 
$x_j$ and $(1 + \beta x_j)$
with $j \neq i$ or $i+1$.
Next, we analyze these polynomials. 
Let us start with $g^T$:
\begin{lemma}
$g^T$ is invariant within the family.
Moreover, $s_ig^T = g^T$.
\end{lemma}

\begin{proof}
Clearly, changing free entries
will not affect powers of $x_j$ and $(1 +\beta x_j)$
with $j \neq i$ or $i+1$.
Let us focus on powers of 
$x_i$, $x_{i+1}$, $(1 + \beta x_i)$ 
and $(1 + \beta x_{i+1})$.
Each ignorable pair contributes
$x_ix_{i+1}$.
Now, consider a frozen $i$.
The column on its right must have
an ignorable pair or a frozen $i+1$.
In either case, 
if we look at the entry the entry above it
and the entry on its top right:
$$
\begin{ytableau}
a & b\\
i & \none
\end{ytableau}
$$
We must have $a > i + 1 \geq b$.
Thus, a frozen $i$ always contributes 
$x_i(1 + \beta x_{i+1}) = X_i$.
Similarly, a frozen $i+1$ 
always contributes
$X_{i+1}$.
Thus, each frozen pair
contributes $X_iX_{i+1}$.

Now, we still need to look at
contributions of $(1 + \beta x_i)$ and 
$(1 + \beta x_{i+1})$ by entries
that are not $i$ or $i+1$.
Assume $j$ is an entry that contributes
$(1 + \beta x_{i+1})$ 
and $j$ is not $i$ or $i+1$.
Let $j'$ be the entry above $j$.
Then $j < i$ and $j' > i+1$.
There is a $k'$ on the row of $j'$
such that $k'$ contributes $(1 +\beta x_i)$.
Also, $k'$ is weakly left of $j'$.
The diagram looks like:
$$
\begin{ytableau}
k' & \dots & j'\\
k  & \dots & j \\
\end{ytableau}
$$
with $k' > i + 1$ and $k < i$.
We pair this $j$ with $k'$.
Similarly, 
given such $k'$,
we can find its corresponding $j$.
In other words,
we pair $(1 + \beta x_i)$ contributors
with $(1 + \beta x_{i+1})$ contributors
that are not $i$ or $i+1$.
This pairing is clearly invariant
under changing free entries. 
\end{proof}

Due to this result, 
we may change our notation $g^T$
into $g^\mathcal{F}$, 
since it only depends on $\mathcal{F}$.
The next step is to study each $f_j^T$.
Clearly, a free $i$ contributes
either $x_i$ or $X_i$.
How can we determine its contribution?
Consider the following lemma:
\ytableausetup{boxsize=10pt,aligntableaux=bottom}

\begin{lemma}
\label{free_i_contribution}
Choose a free $i$ in $T$.
If it is not the last entry
in its block,
then it contributes $x_i$ iff
it is contiguous to the next free $i$.
If it is the last entry in its block,
then it contributes $x_i$ iff
one of the following happens:
\begin{enumerate}
\item It is in the highest row.
\item There is a $b$ on its top right:
\begin{ytableau}
\none & b\\
i
\end{ytableau}
with $b > i + 1$.
\end{enumerate}
\end{lemma}

\begin{proof}
First, assume $i$ is not the last
entry in its block.
We study the entry on its right:
\begin{enumerate}
\item The column on its right 
has ignorable pair.
Then we look at
\begin{ytableau}
a & b\\
\textcolor{red}{i} & \none
\end{ytableau}
where our chosen $i$ is red.
We must have $a > i + 1 \geq b$.
This $i$ contributes $X_i$.
\item The column on its right 
has a free $i$
and this free $i$ is in the same row.
Then we have 
\begin{ytableau}
\none & a\\
\textcolor{red}{i} & i
\end{ytableau}
with $a > i+1$,
or our chosen $i$ is in the top row.
In either case, it contributes $x_i$.
\item The column on its right 
has a free $i$
and this free $i$ is not on the same row
as our chosen $i$.
Then we have
\begin{ytableau}
a & b\\
\textcolor{red}{i} 
\end{ytableau}
with $a > i + 1$ and $b \leq i$.
Our chosen $i$ contributes $X_i$.
\end{enumerate}

Now assume $i$ is the last entry 
in its block.
If it is in the top row, 
then it clearly contributes $x_i$.
Otherwise, 
we look at:
\begin{ytableau}
a & b\\
\textcolor{red}{i}
\end{ytableau}
We know $a > i+1$.
If $b$ exists and $b > i+1$,
then clearly our $i$ contributes $x_i$.
Otherwise, our $i$ contributes $X_i$.
\end{proof}
\ytableausetup{boxsize=normal}

Similarly, for $i+1$, we have:
\begin{lemma}
Choose a free $i + 1$ in $T$.
If it is not the first entry
in its block,
then it contributes $x_{i+1}$ iff
it is contiguous to the previous 
free $i + 1$.
If it is the first entry in its block,
then it contributes $x_{i+1}$ iff
there is an $a$ on its lower left 
with $a < i$.
\begin{align*}
\ytableaushort{\none{i\!\!+\!\!1},a}
\end{align*}
\end{lemma}

\ytableausetup{aligntableaux=center}

We omit the proof since it is basically
the same as the previous one.

Now we understand how the free entries 
contribute.
Clearly, the contribution of one block
is independent from other blocks. 
This implication allows us to 
simplify 
$\sum_{T \in \mathcal{F}} x^{WT(T)}$.
In this family $\mathcal{F}$,
there are $a_j + 1$ ways to fill
the block $j$,
where $a_j$ is the number of 
entries in block $j$.
Let $f_j^l$ be the contribution
of this block when the number of 
$(i+1)$'s is $l$.
$l$ ranges between 0 and $a_j$.
Then we have the following:
$$
\sum_{T \in \mathcal{F}} x^{WT(T)}
= g^\mathcal{F}
\prod_{j = 1}^m\left(\sum_{l = 0}^{a_j}f_j^l\right)
$$

Then we have

\begin{lemma}
$$
\sum_{l = 0}^{a_j} f_j^l 
= \pi_i(f_j^0) \textrm{ or }
\pi_i^K(f_j^0)
$$

Moreover, 
take any $T \in \mathcal{F}$ such
that its $j^{th}$ block has an $i+1$.
Then we are in the second case
iff the first $i+1$ in the $j^{th}$ block
of $T$ contributes $X_{i+1}$.
\end{lemma}

\begin{proof}
First, assume block $j$
only has $i$.
Let $u_p$ be the contribution
of the $p^{th}$ free entry.
Then $f_j^0 = u_1 \dots u_{a_j}$
and each $u_p = X_i$ or $x_i$.

We change the first
free $i$ into $i+1$.
By Lemma \ref{free_i_contribution}, 
this change only affects the first
entry's contribution.
Then $f_j^1 = vu_2 \dots u_{a_j}$
with $v = x_{i+1}$ or $X_{i+1}$.
If $a_j = 1$, we are done by Corollary \ref{pi_comb}.
Otherwise, we change the second 
free $i$ into $i+1$.
The second $i+1$ contributes $x_{i+1}$
iff it is contiguous to the first
free entry.
Also, $u_1 = x_i$ iff the first
entry is contiguous to the second entry.
Thus, we know the second entry
contributes $s_i u_1$.
$f_j^2 = v s_i(u_1) u_3 \dots u_{a_j}$.
Continuing this argument, 
we have $f_j^l =
vs_i(u_1 \dots u_{l-1}) u_{l+1} \dots u_{a_j}$.
The proof is finished by invoking
Corollary \ref{pi_comb}.
\end{proof}

By this result,
$\sum_{l = 0}^{a_j}f_j^l$ must be
symmetric in $i$ and $i+1$.
Recall that we have shown $g^\mathcal{F}$
is symmetric in $i$ and $i+1$.
Thus, $\sum_{T \in \mathcal{F}} x^{WT(T)}$
is symmetric in $i$ and $i+1$.
Finally, we have enough results
to prove Theorem \ref{reverse_SSYT_rule}.
\begin{proof}
Let $\alpha$ be a weak composition
with $\alpha_i > \alpha_{i+1}$.
Let $A \df \{ T \in \mathcal{F}:
K_-(T) \leq \alpha\}$
and $B \df \{ T \in \mathcal{F}:
K_-(T) \leq s_i\alpha\}$.

We only need to show 
\begin{align}
\label{main_equation}
\pi_i^K
\left(\sum_{T \in A} x^{WT(T)}\right)
= \sum_{T \in B} x^{WT(T)}
\end{align}

This is clearly true when 
$B = \emptyset$.
Now assume $B \neq \emptyset$. 
If $A = B$,
then $A = B = \mathcal{F}$,
(\ref{main_equation}) is true since
$\sum_{T \in \mathcal{F}} x^{WT(T)}$
is symmetric in $x_i$ and $x_{i+1}$.

Finally, assume $A$ is a proper subset
of $B$.
We can find $\gamma$ with 
$\gamma_i > \gamma_{i+1}$ such that 
elements in $A$ has left
key $\gamma$ and 
elements in $B$ has left
key $s_i\gamma$.
Then $s_i\gamma \leq s_i\alpha$
and $\gamma \leq \alpha$.
By Lemma \ref{family_and_left_key},
$A$ has elements
whose first block only has $i$.
We have:
$$
\sum_{T \in A} x^{WT(T)}   
=\left( g^\mathcal{F} 
\prod_{j = 2}^m\left(\sum_{l = 0}^{a_j}f_j^l\right)\right)
f_1^0
$$

Take $T \in B$.
Consider its $i+1$ in the first block.
There are two possibilities:
It is in the first column, 
or the column on its left has an
ignorable pair.
In either case,
this $i+1$ contributes $X_{i+1}$,
so
$$
\pi_i^K(f_1^0) = \sum_{l = 0}^{a_1}f_1^l
$$

Finally, letting $f = g^\mathcal{F} 
\prod_{j = 2}^m
(\sum_{l = 0}^{a_j}f_j^l)$,
we have 
\begin{equation*}
\begin{split}
\pi_i^K\left(\sum_{T \in A} x^{WT(T)}\right)
= & \pi_i^K(f \: f_1^0)\\
= & f \: \pi_i^K(f_1^0)\\
= & f \sum_{l = 0}^{a_1}f_1^l\\
= & \sum_{T \in B} x^{WT(T)}
\end{split}
\end{equation*}
\end{proof}

\section{Compatible Word rule 
for Lascoux Polynomials}
\label{S:compatible}

In this section we give another rule for Lascoux polynomials involving compatible pairs.
Recall the set $\C$ from \S \ref{SS:FK} and $\T$ from \S \ref{SS:SV recording}.
We would like to focus on smaller subsets
of them.
\begin{defn}
Let $P$ be a decreasing tableau.
Let $\C_P$ be the set 
consisting of all $(a,i) \in \C$ such that:
\begin{enumerate}
\item $a_j \geq i_j$ for all $j$, and
\item When we insert $a$ into an empty
decreasing tableau using Hecke insertion,
we get $P$.
\end{enumerate}
Correspondingly, 
we define $\T_P$ to be the set consisting 
of all $(P,Q) \in \T$ such that
$K_-(L(Q)) \leq K_+(P)$.
\end{defn}

Then we can introduce 
our main result of this section:
\begin{thm}
\label{restricting_Insert}
The restriction of $\Insert$ to $\C_P$
and the restriction of $\RevInsert$ to $\T_P$, give inverse bijections 
between $\C_P$ and $\T_P$.
\end{thm}

Using this result, we have:
\begin{thm}
\label{compatible_word_rule}
$$
\fLb_{K_+(P)} 
= \sum_{(a,i)\in C_P}
\beta^{\extra(a)} x^{\wt(i)}
$$
\end{thm}
\begin{proof}
By Theorem \ref{set-valued_reverse_SSYT_rule},
we have
$$
\fLb_{K_+(P)} = \sum_{(P,Q)\in \T_P}
\beta^{\extra(Q)}x^{\wt(Q)}
$$
Then the proof is finished by applying $\RevInsert$
on the summands and invoking Theorem 
\ref{restricting_Insert}.
\end{proof}

The rest of this section aims to 
prove Theorem \ref{restricting_Insert}. 
More specifically, 
we only need to show
$\Insert(\C_P) \subseteq \T_P$
and $\RevInsert(\T_P) \subseteq \C_P$.

\subsection{Right key of decreasing tableau}
\label{SS:right key}
In this subsection, 
we investigate the right key
of a decreasing tableau.
First, we introduce an efficient
algorithm that computes the right key.

\subsubsection{Right key via $\star$-action}
We start with a definition:

\begin{defn}
Let $\star$ denote the following right action of the monoid of words with letters in the 
set $\Z_{>0}$, on the set of subsets of $\Z_{>0}$.

Let $S \subseteq \mathbb{Z}_{>0}$
and let $m\in \Z_{>0}$. 
Let $m'$ be the smallest number in $S$ of value
at least $m$.
If $m'$ does not exist, 
we let $S \star m = S \sqcup \{m \}$.
Otherwise,
we define $S \star m = (S - \{ m'\})
\sqcup \{ m \}$.

More generally, 
if $w = w_1 \dots w_n$ is a word of
positive integers, 
we define $S \star w =
(S \star w_1) \star (w_2 \dots w_n)$,
and $S \star w = S$ 
if $w$ is the empty word.
\end{defn}

\begin{ex}
We have: 
\begin{equation*}
\begin{split}
\emptyset \star 3414 & = \{1,4\}\\
\{3,4,7\} \star 3414 & = \{1,4,7\}\\
\{3,4,7\} \star 3141 & = \{1,4,7\}\\
\{3,4,7\} \star 3411 & = \{1,4,7\}\\
\end{split}
\end{equation*}
\end{ex}

We will use this action to 
introduce our right key algorithm.
Before that, 
we need to define a relation on words called K-Knuth equivalence,
which is first introduced in \cite{BS}.
\begin{defn}
The K-Knuth relations are:
\begin{align}
\label{E:K Knuth idempotent}
u\:aa\:v & \Kequiv u\:a\:v \\
\label{E:K Knuth braid}
u\:aba\:v & \Kequiv u\:bab\:v \\
\label{E:Knuth left witness}
u\:bac\:v & \Kequiv u\:bca\:v  
\:\:\: ( a < b < c)\\
\label{E:Knuth right witness}
u\:acb\:v & \Kequiv u\:cab\:v  
\:\:\: ( a < b < c)
\end{align}
where $u, v$ are words and $a<b<c$ are 
positive numbers.

The K-Knuth equivalence $\Kequiv$ is
the transitive and symmetric closure
of these four relations. 

The Knuth equivalence relation $\equiv$ differs from $\Kequiv$ by
removing the relation \eqref{E:K Knuth idempotent} and replacing the 
braid-like relation \eqref{E:K Knuth braid} by the relations
\begin{align}
\label{E:Knuth equal left witness}
 u\:bba\:v & \equiv u\:bab\:v \\
 \label{E:Knuth equal right witness}
 u\:baa\:v & \equiv u\:aba\:v.
\end{align}
\end{defn}

\begin{rem}\label{R:Knuth and star}
It follows from the definitions that $S\star w$ may be obtained
by taking the single column 
semistandard tableau defined by $S$, applying the usual Schensted 
column insertion of $w_1$, then $w_2$, up to $w_n$, and keeping only the first 
column. In particular, applying the well-known fact \cite{Kn} that Knuth-equivalent words 
have the same Schensted insertion tableau, letting $\rev(w)$ be the reverse of the word $w$, 
if $\rev(w)\equiv \rev(w')$ then $S\star w = S\star w'$.
\end{rem}

Words in the same K-Knuth
class have the same $\star$ action.
\begin{lemma}
\label{star_and_K-Knuth}
$S \star w = S \star w'$ if $w \Kequiv w'$.
\end{lemma}
\begin{proof} We may assume that $u$ and $v$ are the empty word.
The result is obvious when $w=aa$ and $w'=a$. If $w$ and $w'$ are related by \eqref{E:Knuth left witness} or
\eqref{E:Knuth right witness} then \eqref{star_and_K-Knuth} holds by Remark \ref{R:Knuth and star}.
For \eqref{E:K Knuth braid} we have
\begin{align*}
S * aba = S * aab = S * ab = S* abb = S* bab
\end{align*}
where the middle two equalities hold by the case of \eqref{E:K Knuth idempotent} and 
the first and last hold by the reverses of the Knuth relations \eqref{E:Knuth equal right witness} and 
\eqref{E:Knuth equal left witness}.
\end{proof}

The $\star$ product is monotonic under left multiplication.

\begin{lem}\label{L:star left multiplication}
For any words $w$ and $v$, $\emptyset \star v \subseteq \emptyset \star wv$ as subsets.
\end{lem}
\begin{proof} We perform an induction on the length of
$v$. When $v$ has length 0, the claim is trivial.

Now assume $v = v'x$ where $x \in \mathbb{Z}_{>0}$.
By the induction hypothesis,
$\emptyset \star v' \subseteq \emptyset \star wv'$.
Next we consider the action of $x$.
If $\emptyset \star v'x\not\subseteq \emptyset \star wv'x$
the only possibility is that 
$y$ in the latter set is replaced by $x$,
but $y$ in the former set is not. 
If this happens, 
$x$ must replace a number smaller than $y$
in the former set, say $z$.
However, $z$ is also in the latter set
and $x$ cannot replace $y$.
Contradiction.
Thus, $\emptyset \star v'x \subseteq \emptyset \star wv'x$.
\end{proof}

We use the $\star$ action to define the right key of a RSSYT and in particular a decreasing tableau.

\begin{defn}\label{D:right key of RSSYT}
For a RSSYT $P$ of partition shape, we define its right key $K_+(P)$ to be the 
RSSYT whose $j$-th column is the column given by 
$\emptyset \star \word(P_{\ge j})$ where $P_{\ge j}$ is the 
decreasing tableau obtained by removing the first $j-1$ columns of $P$.
\end{defn}

\begin{rem} \label{R:right key is a key}
By Lemma \ref{L:star left multiplication} $K_+(P)$ is a key.
\end{rem}

\subsubsection{Right key via Kjdt}
One may also define the right key of a decreasing tableau using Kjdt. 
This is the decreasing analogue of the definition of left key of increasing tableau given in 
\cite{ReY}. We prove the implicit suggestion in \cite{ReY} that the rectification order is irrelevant.

\begin{prop}\label{P:Kjdt right key}
\begin{enumerate}
\item \label{it:1}
For any decreasing tableau $T$ of partition shape, the rightmost column of the 
Kjdt anti-rectification of $T$ with respect to an \textbf{arbitrary} rectification order,
is equal to $\emptyset \star \word(T)$. In particular this column does not depend on the rectification order.
\item \label{it:2}
For any decreasing tableau $P$, the right key of $P$ is the key tableau whose $j$-th column 
is the rightmost column of any Kjdt anti-rectification of $P_{\ge j}$.
\end{enumerate}
\end{prop}
\begin{proof} We only prove part \eqref{it:1} as it immediately implies part \eqref{it:2}.
Let $T'$ be any Kjdt anti-rectification of $T$.
By Theorem 6.2 of \cite{BS},
we know $\word(T') \Kequiv \word(T)$. By Lemma \ref{star_and_K-Knuth},
$\emptyset \star \word(T') = \emptyset\star\word(T)$.
Since $T'$ is a decreasing tableau of antinormal shape, 
$\emptyset \star \word(T')$ agrees with the rightmost column of $T'$.
\end{proof}

\subsubsection{Right key and Hecke insertion}
We determine the precise change in the right key of a decreasing tableau under the operation of 
Hecke insertion of a single value.

We know the following lemma from Theorem 6.2 of \cite{BS}.
\begin{lemma}
Let $P'$ be a decreasing tableau, $x$ a value, and 
$P=(x \HD P')$. 

Then $x\,\word(P') \Kequiv \word(P)$.
\end{lemma}

Then we have:
\begin{lemma}
\label{right_key_change}
Let $P'$ be a decreasing tableau.
Say the insertion $P = (x \HD P')$
ends at column $c$.
Then $K_+(P)$ and $K_+(P')$ agree everywhere
except at column $c$.
Moreover, 
if the insertion causes a contraction,
then $K_+(P) = K_+(P')$.
\end{lemma}
\begin{proof}
Clearly, $P$ and $P'$ agree on column $k$
if $k > c$.
Thus, $K_+(P)$ and $K_+(P')$ agree 
on column $k$.

Now assume $k \leq c$.
Let $w$ (resp. $w'$) be the column word
of $P$ (resp. $P'$) starting at the bottom
of column $k$.
Thus, either $w = w'$ or 
$w \Kequiv yw'$ for some number $y$,
where $y$ is the number inserted to column $k$.
The former case directly implies $K_+(P)$ and $K_+(P')$ 
agree on column $k$.
Now we consider the latter case.
Column $k$ of $K_+(P)$ (resp. $K_+(P')$)
consists of numbers $ \emptyset \star w$ 
(resp. $\emptyset \star w'$).
If $k < c$,
then we know the first character in
$w'$ is less than $y$.
Thus, $\emptyset \star w = \emptyset \star yw' = \emptyset \star w'$.
Finally, 
assume $k = c$ and a contraction occurs.
Then clearly $w = w'$. 
Thus, $K_+(P)$ and $K_+(P')$ agree
on column $c$. 
\end{proof}

What can we say about the changes
at column $c$ of $K_+(P')$ 
if no contraction occurs?
This is answered by the following lemma:
\begin{lemma}
\label{right_key_column_c_change}
Keep the notation from the previous lemma
and assume contraction does not occur.
Let $C$ and $C'$ be column $c$ of $K_+(P)$
and $K_+(P')$ respectively.
Let $D$ be column $c+1$ of $K_+(P)$.
Then as sets $C = C' \sqcup \{e\}$
where $e$ is the smallest number in $C$
that is not in $D$.
\end{lemma}
\begin{ex}
In the following examples we do not 
distinguish between a column and its underlying set.
\begin{enumerate}
\item If $C=\{1,3,4,6,7\}$ and $D=\{1,3,7\}$
then $C'=\{ 1,3,6,7\}$.
\item If $C=\{1,3,4,6,7\}$
and $D$ is empty then $C'=\{ 3,4,6,7\}$.
\end{enumerate}
\end{ex}

\begin{proof}
We prove the lemma by induction on the number
of entries to the right of $C$.
If there are no such entries,
then column $c$ is the rightmost column.
Thus, $C$ (resp. $C'$) agrees with 
column $c$ of $P$(resp. $P'$).
$C'$ has all but the smallest number in $C$,
so the claim holds.

Now let's assume there are
entries to the right of $C$.
Let $m$ be the last number in $word(P)$.
We may pretend as if $m$ does not exist
in $P$ and $P'$ and compute
$C, C'$ and $D$.
By the inductive hypothesis, 
$C' = \{s_1 <  \dots < s_n \}$
and $C = \{s_1 < \dots < s_i 
< e < s_{i+1} < \dots <  s_n\}$,
where $e$ is the extra number.
Moreover, $s_1, \dots, s_{i}$ 
are the $i$ smallest numbers in $D$.
Now we consider the effect of $m$.
There are two cases:
\begin{enumerate}
\item This ignored $m$ is in column $c+1$.
Then column $c+1$ is the last
column of $P$.
Now we let $m$ act on $C$, $C'$ and $D$.
$m$ simply adds itself to $D$.
If $m$ changes $s_j$ into $m$ in $C$,
then $s_j$ is also changed into 
$m$ in $C'$.
Our claim clearly holds.
Otherwise, $e$ in $C$ is changed into $m$.
Then $s_{i+1}$ in $C'$ is changed 
to $m$.
Then $s_{i+1}$ becomes the ``extra number".
Numbers less than it in $C$ are 
$s_1, \dots, s_i$ and $m$.
They are all in $D$.

\item This ignored $m$ is not 
in column $c+1$.
Now we let $m$ act on $C$, $C'$ and $D$.
Assume $m$ changes $s_j$ into $m$ in $C$ with 
$j > i$.
Then $s_j$ is also changed into 
$m$ in $C'$.
For $D$, $m$ will not change $s_1, \dots, s_i$.
Our claim still holds.
Now if $m$ changes $s_j$ 
into $m$ in $C$ with $j \leq i$,
then $m$ also changes $s_j$ in $C'$ 
and $D$ into $m$.
Our claim still holds.
Finally, assume $m$ changes $e$ 
into $m$ in $C$.
Then $m$ will change $s_{i+1}$ in $C'$.
For $D$, $m$ will change a number other than
$s_1, \dots, s_i$.
Then $s_{i+1}$ becomes the ``extra number".
Numbers less than it in $C$ are 
$s_1, \dots, s_i$ and $m$.
They are all in $D$.
\end{enumerate}
\end{proof}

This lemma leads to the following. 
\begin{cor}
Entries in $C'$ are entrywise
weakly less than corresponding
entries in $C$.
\end{cor}

\subsection{Left key of RSSYT}
In this subsection, 
we are going to derive some 
results about the left key.
Our results are analogous 
to the results in the previous subsection.
We start with a result similar to 
Lemma \ref{right_key_change}.

\begin{lemma}
\label{left_key_change}
Take $(a,i) \in \C$.
Assume $a = a'x$ and $i = i'y$,
where $x, y \in \mathbb{Z}_{>0}$.
Let $(P,Q) = \Insert((a,i))$
and $(P',Q') = \Insert((a',i'))$.
Assume the insertion of $x$ ends at column $c$.
Then $K_-(L(Q))$ and $K_-(L(Q'))$ agree everywhere
except at column $c$.
If the insertion of $x$ causes a contraction,
then $K_-(L(Q)) = K_-(L(Q'))$.
\end{lemma}
\begin{proof}
When a contraction occurs, 
clearly $L(Q) = L(Q')$ 
and the conclusion is trivial.
Now assume there's no contraction.
Then $L(Q)$ is obtained by appending 
a number $y$ beneath column $c$ of $L(Q')$.
Then clearly the first $c-1$ columns of 
$K_-(L(Q))$ and $K_-(L(Q'))$ agree.
Consider column $c'$ where $c' > c$.
We know any number 
in column $c'$ of $L(Q)$
is strictly larger than $y$.
Thus, when we compute 
column $c'$ of $K_-(L(Q))$,
this $y$ will be ignored. 
\end{proof}

As in the previous subsection, 
we would like to know 
what happens at column $c$ of $K_-(L(Q'))$
and $K_-(L(Q))$ when a contraction 
does not occur.
The following lemma is an analogue
of \ref{right_key_column_c_change}:

\begin{lemma}
\label{left_key_column_c_change}
Keep the notation from the previous lemma
and assume contraction does not occur.
Let $C$ and $C'$ be column $c$ of 
$K_-(L(Q))$ and $K_-(L(Q'))$ respectively.
Then  $C = C' \sqcup \{e\}$ where for $c=1$, $e=y$ and for 
$c\ge2$, $e$ is the smallest number in $D$ that is not in $C'$
where $D$ is column $c-1$ of $K_-(L(Q))$.
\end{lemma}
\begin{ex}
We have the following examples:
\begin{enumerate}
\item If $C'=\{1,3,6,7\}$
and $D=\{1,3,4,6,7,8\}$ then $C=\{ 1,3,4,6,7\}$.
\item If $C'$ is empty
and $D=\{1,3,4,6,7,8\}$ then $C=\{ 1\}$.
\end{enumerate}
\end{ex}

\begin{proof}
If $c = 1$, our claim is immediate.
Thus, we assume $c \geq 2$.
Let $C_1, C_2, \dots$ be the columns of $L(Q)$ and let
$C_1', C_2', \dots$ be the columns of $L(Q')$.
Then $C = C_1 \triangleleft \dots \triangleleft C_c$,
$C' = C_1' \triangleleft \dots \triangleleft C_c'$
and $D = C_1 \triangleleft \dots \triangleleft C_{c-1}$.
We prove our claim by induction on $c$.

For the base case, we assume $c = 2$.
Then $D = C_1$.
Assume $C_2'=\{s_1 < \dots < s_m\}$.
Then $C_2=\{y < s_1 < \dots < s_m\}$.
Assume $C_1=\{t_1 < t_2 < \dotsm\}$.
Then we know $y \leq t_1$.
When we compute $C_1' \triangleleft C_2'$,
consider two cases:
\begin{enumerate}
\item Case 1: $s_j$ picks $t_j$ for all $j \in [m]$.
When we compute $C_1 \triangleleft C_2$,
$y$ picks $t_1$.
$s_1 < s_2 \leq t_2$, 
so $s_1$ picks $t_2$.
Consequently, $s_j$ picks $t_{j+1}$ for all $j \in [m]$.
Thus, $C_1' \triangleleft C_2'$ contains $t_1, \dots t_m$
while $C_1 \triangleleft C_2$ contains $t_1, \dots t_m, t_{m+1}$.
Our claim is immediate.
\item Case 2: Otherwise let $j$ be the smallest such that
$s_j$ in $C_2'$ does not pick $t_j$.
Thus, $C_1' \triangleleft C_2'$
has $t_1, \dots, t_{j-1}$ but 
does not have $t_j$. 
When we compute $C_1 \triangleleft C_2$,
similar to the previous case,
$y, s_1, \dots, s_{j-1}$ will pick
$t_1, \dots, t_j$.
Then $s_j$ will make the same
choice as in $C_2'$.
Since then, numbers in $C_2$ will
make the same choices as in $C_2'$.
Thus, $C_1 \triangleleft C_2$
has all numbers in $C_1' \triangleleft C_2'$
together with $t_j$.
Our claim is proved. 

\end{enumerate}

Now we do the inductive step.
We first ignore $C_1$ and compute $C$, 
$C'$ and $D$.
By the inductive hypothesis, 
we assume $C' = \{s_1 <  \dots < s_n \}$
and $C = \{s_1 < \dots < s_i 
< e < s_{i+1} < \dots <  s_n\}$,
where $e$ is the extra number.
Assume $D = \{t_1 < t_2 < \dotsm\} $.
Then we know $t_j = s_j$ for all $j \leq i$
and $e = t_{i+1}$.
Now we consider the effect of $C_1$.
We need to study $C_1 \triangleleft C$, 
$C_1 \triangleleft C'$, and $C_1 \triangleleft D$.
$s_1, \dots, s_i$ make the same choices
in all 3 scenarios. 
Notice that since $K_-(L(Q))$ is a key,
numbers in $C_1 \triangleleft C$ must appear
in $C_1 \triangleleft D$.
Thus, when we study 
$C_1 \triangleleft C$,
we can ignore numbers 
in $C_1$ not picked by $D$.
The same is true for $C_1 \triangleleft C'$.
Now we study two cases:
\begin{enumerate}
\item $s_j$ in $C'$ and $t_j$ in $D$ 
make the same choices for all $i < j \leq m$.
In $C_1 \triangleleft C$,
$e$ picks what $t_1$ picks.
Consequently, $s_j$ picks what $t_{j+1}$ picks 
for all $i < j \leq m$.
Our claim is immediate.
\item Case 2: Otherwise
let $j$ be the smallest such that
$s_j$ in $C'$ does not pick what $t_j$ picks.
Thus, $s_{i+1}, \dots, s_{j-1}$ in $C'$ make the same
choices as $t_{i+1}, \dots, t_{j-1}$.
When we compute $C_1 \triangleleft C$,
similar to the previous case,
$e, s_{i+1}, \dots, s_{j-1}$ will make the same
choices as $t_{i+1}, \dots, t_j$.
Then $s_j$ will make the same
choice as in $C'$.
Since then, numbers in $C$ will
make the same choices as in $C'$.
Our claim is proved. 
\end{enumerate}

\end{proof}

\begin{cor}
Values in $C'$ are less than or equal to
corresponding values in $C$.
\end{cor}

\subsection{Proof of Theorem \ref{compatible_word_rule}}

Now we can prove Theorem \ref{compatible_word_rule}.
It is enough to prove the following:
\begin{lemma}
\label{compatible_word_mainlemma}
\begin{enumerate}
\item Take $(a,i) \in \C$.
Assume $a = a'x$ and $i = i'y$,
where $x, y \in \mathbb{Z}_{>0}$.
Let $(P,Q) = \Insert((a,i))$
and $(P',Q') = \Insert((a',i'))$.
Assume $K_+(P') \geq K_-(L(Q'))$.
Assume the insertion of $x$
ends at column $c$ 
and does not cause a contraction.
Then column $c$ of $K_+(P)$ is entrywise
greater than or equal to
column $c$ of $K_-(L(Q))$.
\item Take $(P,Q) \in \T_P$.
Find the smallest number in $Q$
and break ties by picking the 
rightmost.
Suppose it is $y$ in column $c$, 
living in an entry that only contains it.
We remove that entry in $Q$ 
and invoke Hecke reverse insertion on the 
corresponding entry of $P$ with $\alpha = 1$.
Let $x$ be the output.
Then we must have $x \geq y$.
Moreover, assume we get $(P', Q')$ after
the process. 
Then column $c$ of $K_+(P')$ is entrywise
greater than or equal to
the column $c$ of $K_-(L(Q'))$.
\end{enumerate}
\end{lemma}

Why these two statements are enough?
Well, we can use induction to to prove 
$\Insert(\C_P) \subseteq \T_P$
and $\RevInsert(\T_P) \subseteq \C_P$.
For the former, 
we keep all notations in the 
first part of Lemma 
\ref{compatible_word_mainlemma}.
Then clearly $(a',i') \in \C_{P'}$.
By induction on length of $a$,
we may assume $K_+(P') \geq K_-(L(Q'))$
and need to show $K_+(P) \geq K_-(L(Q))$.
The only place where things can go wrong 
is at column $c$ with no contraction occurs. 
Thus, studying this case is enough.
The latter is similar. 
Part 2 of Lemma \ref{compatible_word_mainlemma}
guarantees $(P', Q') \in \T_{P'}$.
Then by inductive hypothesis,
$\RevInsert(P', Q') \in \C_{P'}$. 
Appending $x, y$ respectively makes
the pair in $\C_P$.

Now, we prove Lemma 
\ref{compatible_word_mainlemma}:
\begin{proof}
We begin with the first part.
Keep all notation and assumptions
in the first part of the lemma.
We proceed by considering two cases.
First assume $c = 1$.
Then column 1 of $K_-(L(Q))$ 
(resp. $K_-(L(Q'))$) agrees with
column 1 of $L(Q)$ (resp. $L(Q')$).
Thus, column 1 of $K_-(L(Q))$
is obtained by appending $y$ on the bottom
of column 1 of $K_-(L(Q'))$.
Since column 1 of $K_+(P')$ is entry-wise
less than or equal to 
column 1 of $K_+(P)$,
we only need to worry about 
the new entry on the bottom.
Clearly, the bottom entry in column 1
of $K_+(P)$ equals a number in $P$, 
which is a number in $a$.
It cannot be smaller than $y$, 
which is the smallest number in $i$.
Thus, we are done.

Now assume $c > 1$.
Let $t_1 < \dots < t_m$ be the numbers
in column $c$ of $K_+(P')$.
Let $s_1 < \dots < s_m$ be the numbers
in column $c$ of $K_-(L(Q'))$.
Then we know $t_j \geq s_j$.
Let $e$ be the extra number in 
column $c$ of $K_-(L(Q))$
and assume it is the $i^{th}$ 
smallest number in this column.
By Lemma \ref{left_key_column_c_change},
if $j \leq i$, 
we have the following:
\begin{equation}
\begin{split}
& \: j^{th} 
\textrm{ smallest number in column }
c \textrm{ of } K_-(L(Q)) \\
= & \: j^{th} 
\textrm{ smallest number in column }
c - 1 \textrm{ of } K_-(L(Q)) \\
\leq & \: j^{th} 
\textrm{ smallest number in column }
c - 1 \textrm{ of } K_+(P)\\
\leq & \: j^{th} 
\textrm{ smallest number in column }
c \textrm{ of } K_+(P)\\
\end{split}
\end{equation}

Now if $j > i$,
the $j^{th}$ smallest number 
in column $c$ of $K_-(L(Q))$
is $s_{j-1}$.
The $j^{th}$ smallest number 
in column $c$ of $K_+(P)$
is either at least $t_{j-1}$.
Clearly, our inequality
still holds.

Now we prove the second part.
Keep all notation and assumptions
in the second part of the lemma.
First, $y \leq x$ is immediate:
$y$ is the smallest number in $Q$.
If an entry in $P$ is less than $y$,
than the right key on that entry is
also less than $y$, 
which is a contradiction.

Now we only need to prove 
the bound about keys.
Let $t_1 < \dots < t_m$ be the numbers
in column $c$ of $K_+(P)$.
Let $s_1 < \dots < s_m$ be the numbers
in column $c$ of $K_-(L(Q))$.
Then we know $t_j \geq s_j$.
Assume $t_i$ is not in column $c$
of $K_+(P')$.
Take any $j < i$, we have:

\begin{equation}
\begin{split}
& \: j^{th} 
\textrm{ smallest number in column }
c \textrm{ of } K_+(P') \\
= & \: j^{th} 
\textrm{ smallest number in column }
c + 1 \textrm{ of } K_+(P') \\
\geq & \: j^{th} 
\textrm{ smallest number in column }
c + 1 \textrm{ of } K_-(L(Q'))\\
\geq & \: j^{th} 
\textrm{ smallest number in column }
c \textrm{ of } K_-(L(Q'))\\
\end{split}
\end{equation}

Now, if $j \geq i$,
then the $j^{th}$ smallest number 
in column $c$ of $K_+(P')$
is $t_{j+1}$.
The $j^{th}$ smallest number 
in column $c$ of $K_+(L(Q'))$
is $s_{j+1}$ or $s_j$.
Our equality holds in either case.
\end{proof}

\begin{rem} \label{R:semistandard compatible bounded}
There is a semistandard analogue of Theorem \ref{T:compatible bounded}.
To state this result, we need to modify a few definitions:
\begin{enumerate}
\item In the definition of $\Insert$,  $\RevInsert$ and $\C_P$,  
``Hecke column insertion" is replaced by ``RSK column insertion". 
\item We change $\T$ to be the set of pairs $(P, Q)$,
where $P$ and $Q$ are RSSYTs of the same shape.
\item In the definition of $\T_P$, 
$L(Q)$ is replaced by $Q$ and
the $K_+(P)$ is the classical right key
of a RSSYT.
\end{enumerate}
With these definitions, 
Theorem \ref{restricting_Insert}
is the semistandard analogue. 
It can be proved 
by exactly the same argument. 
\end{rem}

\section{K-theoretic analogue of the
reverse complement map}
\label{S:reverse complement}

In this section
we investigate the map $T \mapsto T^\sharp$ of Proposition \ref{P:decreasing to increasing}.
It is a bijection from decreasing tableaux
to increasing tableaux. 
Recall the map is defined as follows.
Take a decreasing tableau $T$.
Find the smallest rectangle
that contains $T$.
Then we do several ``iterations" of Kjdt 
to anti-rectify $T$.
In each iteration, 
we perform Kjdt at the leftmost empty
space of each row, from top to bottom.
Then we rotate the result by 180
degrees and obtain $T^\sharp$. 

\begin{ex}
\label{E:anti-rectify decreasing}
We recompute $\JSE(T)$ from Example \ref{X:Krect} using the above iterations. We start with the decreasing tableau $T$, along with the rectification order (a tableau giving the order that boxes are occupied during the computation 
of $\JSE(T)$), and another tableau whose entry $i$ means that the given box is occupied during the $i$-th iteration.
\[
T = 
\begin{ytableau}
9 & 8 & 7 & 5 & 3\\
7 & 5 & 4 & 3\\
4 & 2 & 1 \\
3\\
1\\
\end{ytableau}
\quad 
\ytableaushort{\none\none\none\none\none,\none\none\none\none1,\none\none\none25,\none368a,\none479b}
\quad
\ytableaushort{\none\none\none\none\none,\none\none\none\none1,\none\none\none12,\none1234,\none1234}
\]

Let $T^{(i)}$ be the skew tableau just before the $i$-th iteration.
So $T=T^{(1)}$ and $T^{(m)}=\JSE(T)$ where $m$ is the number of columns of $T$.

$T^{(2)}$, $T^{(3)}$, $T^{(4)}$ and 
$T^{(5)}$ are listed below:
$$
\begin{ytableau}
\none & 9 & 8 & 7 & 5\\
\none & 7 & 5 & 4 & 3\\
\none & 4 & 2 & 1 \\
\none & 3\\
9 & 1\\
\end{ytableau}
\:\:\:\:
\begin{ytableau}
\none & \none & 8 & 7 & 5\\
\none & \none & 5 & 4 & 3\\
\none & \none & 4 & 2 & 1 \\
\none & 9 & 3\\
9 & 7 & 1\\
\end{ytableau}
\:\:\:\:
\begin{ytableau}
\none & \none & \none & 7 & 5\\
\none & \none & \none & 5 & 3\\
\none & \none & 8 & 4 & 1 \\
\none & 9 & 5 & 2\\
9 & 7 & 3 & 1\\
\end{ytableau}
\:\:\:\:
\begin{ytableau}
\none & \none & \none & \none & 7\\
\none & \none & \none & \none & 5\\
\none & \none & 8 & 7 & 3 \\
\none & 9 & 5 & 4 & 2\\
9 & 7 & 3 & 2 & 1\\
\end{ytableau}
$$
Once we rotate $T^{(5)}$ by 180 degrees,
we get:
$$
T^\sharp = 
\begin{ytableau}
1 & 2 & 3 & 7 & 9\\
2 & 4 & 5 & 9\\
3 & 7 & 8 \\
5\\
7\\
\end{ytableau}
$$
Consider the rightmost columns of the tableaux $T^{(5)},T^{(4)},\dotsc,T^{(1)}$ in that order. We get
\[
\ytableaushort{75553,5333,311,2,1}
\]
Lemma \ref{skew-right-key-shifts} implies that this is $K_+(T)$.
\end{ex}

Now, we have:
\begin{prop}
\label{reverse_complement_fixes_keys}
Let $T$ be a decreasing tableau.
Then
$$
K_+(T) = K_-(T^\sharp)
$$
\end{prop}
\begin{ex}
In Example \ref{E:anti-rectify decreasing},
$K_+(T)$ and $K_-(T^\sharp)$ are both
$(3,1,5,0,4,0,1)$.
\end{ex} 

The rest of this section aims to prove it.
We start with a few definitions.
\begin{defn}
Let $T$ be a
decreasing skew tableau.
We use $T_{\leq c, \leq r}$
to denote the decreasing skew tableau
obtained by keeping the
first $c$ columns 
and first $r$ rows of $T$. 
Analogously, we define similar notations
such as $T_{\geq c, \leq r}$. 
\end{defn}
\begin{defn}
We generalize $K_+(T)$ to 
decreasing skew tableaux $T$.
Let $K_+(T)$ be the normal shape decreasing tableau whose
$c$-th column equals
$\emptyset \star w(T_{\geq c, < \infty})$ for all $c$.
\end{defn}

\begin{ex}
\label{generalized-right-key-example}
Let $T^{(2)}$, \dots, $T^{(5)}$ 
be the tableaux in 
Example \ref{E:anti-rectify decreasing}.
Their right keys are:
$$
\begin{ytableau}
7 & 7 & 5 & 5 & 5\\
5 & 5 & 3 & 3 & 3\\
3 & 3 & 1 & 1 \\
2 & 2\\
1 & 1\\
\end{ytableau}
\:\:\:\:
\begin{ytableau}
7 & 7 & 7 & 5 & 5 \\
5 & 5 & 5 & 3 & 3 \\
3 & 3 & 3 & 1 & 1 \\
2 & 2 & 2\\
1 & 1 & 1\\
\end{ytableau}
\:\:\:\:
\begin{ytableau}
7 & 7 & 7 & 7 & 5\\
5 & 5 & 5 & 5 & 3\\
3 & 3 & 3 & 3 & 1\\
2 & 2 & 2 & 2\\
1 & 1 & 1 & 1\\
\end{ytableau}
\:\:\:\:
\begin{ytableau}
7 & 7 & 7 & 7 & 7\\
5 & 5 & 5 & 5 & 5\\
3 & 3 & 3 & 3 & 3\\
2 & 2 & 2 & 2 & 2\\
1 & 1 & 1 & 1 & 1\\
\end{ytableau}
$$
\end{ex}

\begin{rem}
We make a few observations about this generalization.

\begin{enumerate}
\item $K_+(T)$ is a key.
\item The rightmost column of $K_+(T)$ and $T$
must agree. 
\item When $T$ is of normal shape,
this definition agrees with our previous definition.
\end{enumerate}
\end{rem}

Let us look at Example \ref{generalized-right-key-example}
carefully. 
$T^{(4)}$ is obtained from $T^{(3)}$ 
by applying Kjdt
on column 4 twice. 
Their right keys are the same
except at column 4.
Moreover, 
column 4 of $T^{(4)}$ agrees with
column 3 of $T^{(3)}$.
This phenomenon is captured in
the next result:

\begin{lemma}
\label{skew_right_key_changes}
Let $T$ be a decreasing skew tableau.
Assume its $c$-th column ends at row $r_1$, 
while its $(c+1)$-th column ends at row $r_2$.
Assume $r_1 > r_2$.
We do Kjdt at the leftmost position
of row $r_2 + 1, \dots, r_1$
and let $T'$ be the result. 
Then $K_+(T)$ and $K_+(T')$ 
agree everywhere except at column $c+1$.
Moreover, column $c+1$ of $K_+(T')$
agrees with column $c$ of $K_+(T)$.
\end{lemma}
\begin{proof}
Notice: $T_{>c+1, < \infty} 
= T'_{>c+1, < \infty}$.
Thus, $K_+(T)_{>c+1, < \infty} 
= K_+(T')_{>c+1, < \infty}$.

Now assume $i \leq c + 1$ and 
column $i$ of $T$ ends at row $j$.
Column $i$ of $K_+(T)$ has numbers 
$\emptyset \star w(T_{\geq i, < \infty})
= \emptyset \star w(T_{\geq i, \leq j})$.
Now we look at $w(T_{\geq i, \leq j})$.
It is also 
$w(T_{= i, \leq j})w(T_{> i, \leq j})$.
Because $T_{\leq i, \leq j}$ is reverse
rectified, 
$\emptyset \star w(T_{\leq i, \leq j}) 
= \emptyset \star w(T_{=i, \leq j})$.
Thus, column $i$ of $K_+(T)$ has numbers:
$\emptyset \star w(T_{\leq i, \leq j})
w(T_{> i, \leq j})
= \emptyset \star w(T_{< \infty, \leq j})$

Next, if $i \leq c$, 
then we know column $i$ of $T'$ also
ends at row $j$.
Similar to above,
column $i$ of $K_+(T')$ has numbers
$\emptyset \star w(T'_{< \infty, \leq j})$.
Notice that $T'_{< \infty, \leq j}$
is obtained from $T_{< \infty, \leq j}$ by
a sequence of Kjdt moves.
Thus, their column words are K-Knuth 
equivalent, so they yield the same
results when acting on $\emptyset$.

Finally, we study column $c+1$.
Column $c+1$ of $K_+(T')$ is
$\emptyset \star w(T'_{< \infty, \leq r_1})$.
Column $c$ of $K_+(T)$ is 
$\emptyset \star w(T_{< \infty, \leq r_1})$.
Clearly $T'_{< \infty, \leq r_1}$
is obtained from $T_{< \infty, \leq r_1}$
by K-jdt moves,
so we are done.
\end{proof}

Then we can describe how the right key changes
during the map $T \mapsto T^\sharp$.
In Example \ref{generalized-right-key-example},
each iteration changes the right key
in the following way:
It copies column $c$ to column $c+1$,
where $c$ goes from 4 to 1.
This pattern holds in general:

\begin{lemma}
\label{skew-right-key-shifts}
Let $T$ be a decreasing tableau.
Let $C_1, \dots, C_m$ be
columns of $K_+(T)$.
Assume we finished $t$ iterations
while computing $T^\sharp$.
Then the right key of the current
skew-shape tableau has columns:
$C_1, C_1, \dots, C_1, C_2, \dots, C_{m-t}$
\end{lemma}
\begin{proof}
We prove by induction. 
When $t = 0$, the claim is trivial.
Now assume the claim holds
after $t$ iterations.
Assume now column $i$ ends at $r_i$.
During the $(t+1)^{th}$ iteration,
the algorithm performs Kjdt at the
leftmost position of row $r_m + 1, 
r_m + 2, \dots, r_1$.
How do these moves affect the right key?
By Lemma \ref{skew_right_key_changes},
it will copy column $c$ of the right key
to column $c+1$, 
where $c$ goes from $m - 1$ to 1.
Thus, after this iteration,
the right key becomes:
$C_1, C_1, \dots, C_1, C_2, \dots, C_{m-t - 1}$.
\end{proof}

\begin{proof}[Proof of Proposition \ref{reverse_complement_fixes_keys}]
Let $C_1, \dots, C_m$ be the columns of $K_+(T)$.
By Lemma \ref{skew-right-key-shifts},
after $t$ iterations of the computation of $T^\sharp$,
the rightmost column of the resulting skew
tableau (call it $T^{(t+1)}$) is $C_{m - t}$.
By definition we also have $T^{(m)} = \JSE(T)$.

On the other hand consider the computation of the $(m-t)$-th column $C$ of $K_-(T^\sharp)$.
By definition we take the first $m-t$ columns of $T^\sharp$ within a tight rectangle, 
moving to the southeast via a Kjdt which infuses with a column superstandard tableau having $m-t-1$ columns.
$C$ is the leftmost column of the result.
Rotating by 180 degrees, this is the same as taking 
$(T^\sharp)^*=\JSE(T)$, using the last $m-t-1$ columns in a column superstandard Kjdt to move it to the northwest, and 
taking the rightmost column $C$. However the $m-t-1$ iterations which slide into these columns are precisely undoing the last $m-t-1$ iterations in the computation of $\JSE(T)$. Therefore $C$ equals the rightmost column of 
$T^{(t+1)}$, which is $C_{m-t}$, as required.
\end{proof}

\begin{ex}
We illustrate this argument on column 
3 of $K_+(T)$,
where $T$ is the tableau in 
Example \ref{E:anti-rectify decreasing}.
By definition, 
here's how we compute column 3 of $K_-(T^\sharp)$.
We start with $T^\sharp$.
Then we Kjdt to anti rectify $T^\sharp$.
We stop once the first 3 columns are
anti rectified:
$$
\begin{ytableau}
\none & \none & 1 & 7 & 9\\
\none & \none & 3 & 9\\
\textcolor{red}{1} & 2 & 4 \\
\textcolor{red}{3} & 4 & 5\\
\textcolor{red}{5} & 7 & 8\\
\end{ytableau}
$$
Then the red numbers form column 3
of $K_-(T^\sharp)$.
Notice that when we rotate this tableau,
we get $T^{(3)}$ 
in Example \ref{E:anti-rectify decreasing},
whose rightmost column is column $3$ of 
$K_+(T)$.

\end{ex}

\section{Alternative descriptions of $K_+(T)$}
\label{S:alternative right key}
This section provides some alternative 
descriptions of the right key $K_+(T)$,
where $T$ is a decreasing tableau.
\begin{prop}
\label{P: Alternative definitions of the right key}
Let $T$ be a decreasing tableau.
Let $T_{\ge j}$ be the decreasing 
tableau obtained by removing the first $j-1$ columns of $T$.
Then the column $C$ computed 
by the following procedures all agree with
column $j$ of $K_+(T)$.
\begin{enumerate}
\item $C =\emptyset \star \word(T_{\ge j})$.
\item $C$ is the rightmost column of a \textbf{arbitrary} anti-rectification of 
$T_{\ge j}$.
\item Assume column $j$ of $T$ has $m$ numbers.
We conjugate $T_{\ge j}$ and obtain $T_{\ge j}'$.
Then we invoke Hecke reverse column insertion
at the end of column $m, \dots, 1$ with $\alpha = 1$.
Then $C$ is the set of output numbers. 
\end{enumerate}
\end{prop}

\begin{proof}
Let $T$ be a decreasing tableau. 
Proposition \ref{P:Kjdt right key} shows that the first two methods for 
computing $K_+(T)$ agree.
Next, we show procedure 3 agrees with
procedure 1. 
We may assume $j = 1$.
Thus we need to prove:
Assume we conjugate $T$ and get $T'$.
Invoke reverse insertion on each column 
of $T'$ from right to left. 
Show the output agree with $\emptyset \star \word(T)$.

Prove by induction on $m$,
the number of rows in $T$.
When $m = 1$, 
our claim is trivial.
Now we do the inductive step.
Assume both column 1 and column 2 of $T$ 
has $m$ entries. 
Notice that the result of reverse
insertions and $\emptyset \star \word(T)$ are not affected
if we ignore column 1 of $T$.
Thus, by keep removing columns on the left, 
we may assume $T$ has only 
one column with $m$ entries.
Now, we invoke reverse Hecke insertion
at the end of column $m$ of $T'$
and let $x_m$ be the output.
We conjugate the resulting tableau and get 
$P$.
Then $\word(T') \equiv_K x_m \word(P')$.
Define the row word of $T$ to be the following:
Read entries of $T$ from bottom to top,
within each column from left to right.
Then $\rev(\word(T'))$ is the row
word of $T$.
By Theorem 5.4 of \cite{BS}
the row word of $T$ is K-Knuth equivalent
to $\word(T)$.
Thus, 
we know 
$$\emptyset \star \word(T) = \emptyset \star \rev(\word(T')) 
= \emptyset \star \rev(\word(P'))x_m$$
Similarly, $\rev(\word(P'))$ is
the row word of $P$, so
$$
\emptyset \star \word(P) = \emptyset \star \rev(\word(P'))
$$

Next, invoke reverse insertion
at column $m-1, m-2, \dots, 1$ of $P'$
and get numbers $x_{m-1}, \dots, x_1$.
Since $P$ has $m-1$ rows,
our inductive hypothesis says
$$
\emptyset \star \word(P) = \{x_{m-1}, \dots, x_1\}
$$
By Lemma \ref{Pieri property of Hecke insertion},
$x_m > x_{m-1} > \dots > x_1$.
Thus, 
\begin{equation*}
\begin{split}
\emptyset \star \word(T) 
=  & \emptyset \star \rev(\word(P'))x_m \\
=  & \emptyset \star \word(P)x_m \\
=  & \{x_{m-1}, \dots, x_1\} \star x_m \\
=  & \{x_{m}, \dots, x_1\} 
\end{split}
\end{equation*}
\end{proof}


\begin{thebibliography}{LLMSSZ}

\bibitem[AS]{AS} 
S. Assaf and A. Schilling.
A Demazure crystal construction for Schubert polynomials.
Algebr. Comb. 1 (2018) no. 2, 225--247. 

\bibitem[B]{B:Gr} 
A. Buch.
A Littlewood-Richardson rule for the K-theory of Grassmannians. 
Acta Math. 189 (2002), no. 1, 37--78.

\bibitem[B2]{B}
A. Buch.
Alternating signs of quiver coefficients. 
J. Amer. Math. Soc. 18 (2005), no. 1, 217--237. 

\bibitem[BF]{BF} 
A. Buch and W. Fulton. 
Chern class formulas for quiver varieties. 
Invent. Math. 135 (1999), no. 3, 665--687. 

\bibitem[BJS]{BJS}
S. Billey, W. Jockusch, and R. Stanley.
Some combinatorial properties of Schubert polynomials.
J. Algebraic Combin. 2 (1993), no. 4, 345--374. 

\bibitem[BKSTY]{BKSTY} A. Buch, A. Kresch, M. Shimozono, H. Tamvakis, and A. Yong.
Stable Grothendieck polynomials and K-theoretic factor sequences. 
Math. Ann. 340 (2008), no. 2, 359--382.

\bibitem[BS]{BS}
A.S. Buch and M.J. Samuel, 2016. K-theory of minuscule varieties. 
Journal f\"ur die reine und angewandte Mathematik, 2016 (719), pp. 133--171.

\bibitem[BSS]{BSS}
G. Benkart, J. Stroomer, and F. Sottile.
Tableau switching: algorithms and applications.
J. Combin. Theory Ser. A 76 (1996), no. 1, 11--43. 

\bibitem[BSW]{BSW}
V. Buciumas, T. Scrimshaw and K. Weber.
Colored five-vertex models and Lascoux polynomials and atoms.
J. Lond. Math. Soc. (2) 102 (2020), no. 3, 1047--1066. 


\bibitem[Bu]{Bu} A. Buch. 
A Littlewood-Richardson rule for the $K$-theory of 
Grassmannians.
Acta Math. 189 (2002), no. 1, 37--78. 

\bibitem[EG]{EG}
P. Edelman and C. Greene.
Balanced tableaux. Adv. in Math. 63 (1987), no. 1, 42--99.

\bibitem[FK]{FK} 
S. Fomin and A. N. Kirillov.
Grothendieck polynomials and the Yang-Baxter equation. 
Formal power series and algebraic combinatorics/S\'eries formelles et combinatoire alg\'ebrique, 
DIMACS, 1994, 183--189, Piscataway, NJ.

\bibitem[GLS]{GLS} M. Gurevich, E. Lapid, and M. Shimozono.
Robinson-Schensted-Knuth correspondence in the representation theory of the general linear group over a non-archimedean local field. To appear in Represent. Theory.

\bibitem[H]{H} 
M. Haiman.
On mixed insertion, symmetry, and shifted Young tableaux.
J. Combin. Theory Ser. A 50 (1989), no. 2, 196--225. 

\bibitem[Kir]{Kir}
A. N. Kirillov.
Notes on Schubert, Grothendieck and key polynomials. S
SIGMA Symmetry Integrability Geom. Methods Appl. 
12 (2016), Paper No. 034, 56 pp.

\bibitem[KMS]{KMS} 
A. Knutson, E. Miller, and M. Shimozono.
Four positive formulae for type A quiver polynomials.
Invent. Math. 166 (2006), no. 2, 229--325. 

\bibitem[KMY]{KMY} A. Knutson, E. Miller, and A. Yong.
Tableau complexes.
Israel J. Math. 163 (2008), 317--343. 

\bibitem[KMY2]{KMY2} 
A. Knutson, E. Miller, and A. Yong.
Gr\"obner geometry of vertex decompositions and of flagged tableaux. 
J. Reine Angew. Math. 630 (2009), 1--31.

\bibitem[Kn]{Kn} D. E. Knuth. 
Permutations, matrices, and generalized Young tableaux.
Pacific J. Math. 34 (1970) no. 3, 709--727.

\bibitem[LamSh]{LS:QHGr} 
T. Lam and M. Shimozono.
Quantum cohomology of $G/P$ and homology of affine Grassmannian. 
Acta Math. 204 (2010), no. 1, 49--90.

\bibitem[Las]{Las} 
A. Lascoux.
Schubert \& Grothendieck: un bilan bid\'ecennal, S\'em. Lothar. Combin. 50 (2003/04), Art. B50i.

\bibitem[LasSc]{LS:Groth}
A. Lascoux and M.~P.~Sch\"utzenberger. Structure de Hopf de l'anneau de cohomologie et de l’anneau de Grothendieck d'une vari\'et\'e de drapeaux. C. R. Acad. Sci. Paris S\'er. I Math. 295, 11 (1982), 629--633.

\bibitem[LasSc2]{LS:Schub}
A. Lascoux M.~P.~Sch\"utzenberger. 
Tableaux and non-commutative Schubert polynomials.
Funct. Anal. Appl. 23 (1989), 63--64.

\bibitem[LasSc3]{LS:keys}
A. Lascoux M.~P.~Sch\"utzenberger. 
Keys and standard bases.
Invariant theory and tableaux (Minneapolis, MN, 1988), 125--144,
IMA Vol. Math. Appl., 19, Springer, New York, 1990. 

\bibitem[LLMSSZ]{LLMSSZ}
 T.~Lam, Luc~Lapointe, J.~Morse, A.~Schilling, M.~Shimozono, and M.~Zabrocki.
k-Schur functions and affine Schubert calculus. Fields Institute Monographs, 33. 
Springer, New York; Fields Institute for Research in Mathematical Sciences, Toronto, ON, 2014. 

\bibitem[LLS]{LLS}
T. Lam, S. J. Lee, and M. Shimozono. Back stable Schubert calculus.
 Compos. Math. 157 (2021), no. 5, 883--962. 
 
\bibitem[LLS2]{LLS:back stable Groth}
T. Lam, S. J. Lee, and M. Shimozono. Back stable $K$-theory Schubert calculus.
In preparation.

\bibitem[Mat]{Mat} 
T. Matsumura.
Flagged Grothendieck polynomials.
J. Algebraic Combin. 49 (2019), no. 3, 209--228.

\bibitem[Mi]{Mi}
E. Miller, 
Alternating formulas for K-theoretic quiver polynomials. 
Duke Math. J. 128 (2005), no. 1, 1--17.
 
\bibitem[Mon]{Mon}
C. Monical. 
Set-valued skyline fillings. 
S\'em. Lothar. Combin. 78B (2017), Art. 35, 12 pp. 

\bibitem[PeSc]{PS}
O. Pechenik and T. Scrimshaw. 
K-theoretic crystals for set-valued tableaux of rectangular shapes. 
S\'em. Lothar. Combin. 84B (2020), Art. 19, 12 pp. 

\bibitem[Pet]{Pet}
D. Peterson. Quantum cohomology of $G/P$. Lecture notes, MIT, 1997.

\bibitem[ReY]{ReY} V. Reiner and A. Yong. 
The ``Grothendieck-to-Lascoux" conjecture.
preprint arXiv:2102.12399v1, 2021.

\bibitem[RS]{RS}
V. Reiner and M. Shimozono.
Key polynomials and a flagged Littlewood-Richardson rule. 
J. Combin. Theory Ser. A 70 (1995), no. 1, 107--143. 

\bibitem[RoY]{RoY} C. Ross and A. Yong.
Combinatorial rules for three bases of polynomials.
S\'em. Lothar. Combin. 74 (2015), Article B74a, 11pp.

\bibitem[TY]{TY} H. Thomas and A. Yong.
A jeu de taquin theory for increasing tableaux, with applications 
to K-theoretic Schubert calculus. 
Algebra Number Theory 3 (2009), no. 2, 121--148. 

\bibitem[Wil]{Willis} 
M. Willis.
A Direct Way to Find the Right Key of a Semistandard Young Tableau.
Annals of Combinatorics, 17, no. 2, p. 393--400, 2013.




\end{thebibliography}
\end{document}